\documentclass[preprint]{elsarticle}
\usepackage{amssymb}
\usepackage{amsmath}
\usepackage{amssymb,amsmath,mathrsfs,bm,xcolor}
\usepackage{cases}
\usepackage{array}
\usepackage{wasysym,latexsym}
\usepackage{graphicx,color,algorithm,algpseudocode}
\usepackage{booktabs,makecell,multirow}
\usepackage{epstopdf,epsfig}
\usepackage{diagbox}
\usepackage{float}
\usepackage{mathtools}
\usepackage{setspace}
\usepackage{amsthm}

\newtheorem{thm}{Theorem}
\newtheorem{lemma}{Lemma}
\newtheorem{remark}{Remark}

\usepackage{setspace}
\usepackage{graphicx} 

\usepackage{algorithm,algpseudocode}
\usepackage{algpseudocode} 

\usepackage{caption} 
\usepackage{subcaption} 
\usepackage{adjustbox} 
\usepackage[absolute,overlay]{textpos} 
 \usepackage{geometry} 

\usepackage[colorlinks,linkcolor=blue,anchorcolor=blue,linktocpage=true,urlcolor=blue,citecolor=blue]{hyperref}

\usepackage{enumitem} 
\usepackage{hyperref} 

\newcommand{\msc}[1]{\textbf{MSC Code}#1}

\journal{ }
\graphicspath{{figures/}}
\begin{document}
	
	\begin{frontmatter}
		
		
		
 \title{New Primal-Dual Algorithm  for Convex Problems}
 \author[1]{Shuning Liu}
\ead{snliu00@163.com}

\author[2]{Zexian Liu\corref{cor1}}
\ead{liuzexian2008@163.com}

\address[1]{School of Mathematics and Statistics, Guizhou University, Guiyang 550025, China}
\address[2]{School of Mathematics and Statistics, Guizhou University, Guiyang 550025, China}

\cortext[cor1]{Corresponding author}

\begin{abstract}
  Primal-dual algorithm (PDA) is a classic and  popular scheme  for   convex-concave saddle point problems. It is universally acknowledged that  the proximal terms in the subproblems about the primal and dual variables are crucial to   the convergence theory and  numerical performance of primal-dual algorithms. By taking advantage of  the information from the current and previous iterative points, we exploit  two new  proximal terms  for the subproblems about
the primal and dual variables. Based on two new     proximal terms,   we present a new primal-dual algorithm for   convex-concave saddle point problems with bilinear coupling terms  
and establish its   global convergence and  $\mathcal{O}(1/N)$   ergodic convergence rate. When either the primal function or the dual function is strongly convex, we accelerate the above proposed algorithm and  show that the corresponding  algorithm can achieve $\mathcal{O}(1/N^2)$ convergence rate. Since the conditions for the stepsizes of the proposed algorithm are related directly to the spectral norm of the linear transform,   which is difficult to obtain in some applications,   we  also  introduce  a linesearch strategy for the above  proposed  primal-dual algorithm   and establish its global convergence and  $\mathcal{O}(1/N)$ ergodic convergence rate   . 
Some numerical experiments are conducted on matrix game and LASSO problems   by comparing with other state-of-the-art  algorithms, which  demonstrate  the effectiveness of the proposed three primal-dual algorithms.
\end{abstract}		
		\begin{keyword}
        primal-dual algorithms, saddle-point problems,  linesearch, global convergence,  convergence rate			
		\end{keyword}      
	\end{frontmatter}
\msc{: 49M29, 65K10, 65Y20, 90C25}
\section{Introduction.}
    \par Let $X$ and  $Y$ be two finite-dimensional Euclidean spaces with an inner product $\langle \cdot , \cdot \rangle$ and the induced norm  $\| \cdot \|=\sqrt{\langle \cdot , \cdot \rangle}$. In this paper, we focus on the following saddle problem with bilinear  coupling term:
         \begin{equation}\label{p}
            \min_{x\in X} \max_{y\in Y}\;L(x,y)=g(x)+\langle Kx,y\rangle-f^*(y),
        \end{equation}
    where $g:X\rightarrow(-\infty,+\infty]$ and $f:Y\rightarrow(-\infty,+\infty]$ are proper real-valued closed  convex functions,  $f^*$ denotes the conjugate function of  $f$ and is defined by $f^*(y)=\underset{z\in Y}{\sup}\{\langle y,z\rangle-f(z) | y \in Y\}$,   $K:X  \rightarrow  Y $ is a bounded linear operator,   $K^T$ denotes the matrix transpose or adjoint operator of $K$ and   $L=\|K\|$.
    \par Under the above assumptions, the problem (\ref{p}) is equivalent to the primal problem
        \begin{equation}\label{pp1}
            \min_{x\in X}\;f(Kx)+g(x),
        \end{equation}
    and the dual problem
        \begin{equation}\label{pp2}
             \max_{y\in Y}\;f^*(y)+g^*(-K^Tx).
        \end{equation}
    \par Problem (1) has applications in numerous fields, including but not limited to statistics, mechanics, economics, signal processing, image processing, and machine learning, among others. For further details, see \cite{w7,w9,w26,w28}.\\

\subsection{Related algorithms.}
     \par There exist numerous well-known algorithms for solving problems (\ref{p})-(\ref{pp2}) simultaneously, such as primal-dual algorithms \cite{w6,w13,w3,w4} and  their variants \cite{w1,w18,w11,w14,w12,w25}, and the alternating direction method of multipliers (ADMM) \cite{w23,w29,w27,w24} and  their variants \cite{w5}, and some recent algorithms and variants \cite{B1,B2,B3,B4} .
     \par A classical and simple primal-dual full-splitting algorithm \cite{w17} was proposed, which is now  widely known as the Arrow-Hurwicz method  and addresses the saddle problem (\ref{p}) by alternately  minimizing over the primal variable $x$ and maximizing over the dual variable $y$. The Arrow-Hurwicz method is described  as follows: 
         \begin{equation}
             \left\{
             \begin{aligned}
            &x_{k} = \text{Prox}_{\tau g}(x_{k-1} - \tau K^Ty_{k-1}), \\
            &y_{k} = \text{Prox}_{\sigma f^*}(y_{k-1} + \sigma Kx_{k}),
            \end{aligned}
            \right.
        \end{equation}
    where $\tau > 0$ and $\sigma > 0$ are stepsize parameters. However, He et al. \cite{w19} constructed an example to illustrate that the Arrow-Hurwicz method does not always converge.
    \par Chambolle and Pock  \cite{w1, w20}  introduced an extrapolation step $\bar x_k=x_k+\theta(x_{k}-x_{k-1})$ for the  Arrow-Hurwicz method and presented  the primal-dual algorithm (PDA), which  is given by 
        \begin{equation} \label{eq:PDA}
            \left\{
           \begin{aligned}
             &x_{k}=\text{Prox}_{\tau g}(x_{k-1}-\tau K^Ty_{k-1}),\\&
            \bar x_k=x_k+\theta(x_{k}-x_{k-1}),\\&
             y_{k}=\text{Prox}_{\sigma f^*}(y_{k-1}+\sigma K\bar x_k).
        \end{aligned}
        \right.
        \end{equation}
When $\theta = 1$,  the convergence   of PDA was established under the condition $\tau \sigma L^2 < 1$ \cite{w1}. Note that convergence results of PDA with $\theta \in (0, 1)$ issue remains unresolved \cite{w15}.
    \par Based on the golden ratio algorithm \cite{w10} for variational inequalities,  Chang and Yang \cite{w8} proposed a golden ratio primal-dual algorithm (GRPDA), which abandons the extrapolation step and uses a convex combination of   all the previously generated primal iterative points. GRPDA is given by  
         \begin{equation}
           \left\{
            \begin{aligned}
             &z_k=\frac{\psi-1}{\psi}x_{k-1}+\frac{1}{\psi}z_{k-1},\\&
             x_{k}=\text{Prox}_{\tau g}(z_{k}-\tau K^Ty_{k-1}),\\&
             y_{k}=\text{Prox}_{\sigma f^*}(y_{k-1}+\sigma Kx_{k}).
            \end{aligned}
            \right.
         \end{equation}
  The global  convergence and ergodic convergence rate \cite{w8} of GRPDA   were also established    when  $\tau \sigma L^2 < \psi$, where $\psi \in (1,\phi]$ and $\phi=\frac{\sqrt{5}+1}{2}$.
    \par  Malitsky and  Pock \cite{w2} first introduced a linesearch strategy for    PDA  algorithm (PDAL). Similarly, Chang et al. \cite{w15} incorporated a linesearch into GRPDA and presented an efficient GRPDA with linesearch (GRPDAL). 
		
\subsection{Motivations and contributions.}
      \par 

    In PDA   \eqref{eq:PDA}, the subproblem about $x$ is solved by proximal point algorithm (PPA), namely, 
     \begin{align*}
            x_{k+1} =\arg \min_{x \in X} \{g(x)+\langle Kx,y_k\rangle+\frac{1}{2\tau}\|x-x_k\|^2\},
        \end{align*}  
    and so is for the subproblem about $y$.  
    As is well known, the proximal terms $ \frac{1}{2\tau}\left\|x-x_k\right\|^2$ and $-\frac{1}{2\sigma}\left\|y-y_k\right\|^2$ in PPA are of great importance for the convergence analysis  and   numerical performance of PDA. A suitable proximal term can accelerate the corresponding algorithms greatly. One line of research about the proximal terms in PDA  focuses on   designing various efficient strategies for stepsize parameters $\tau$ and $\sigma$. Another line focuses on   developing some new points to replace $x_k$ in proximal term for  speeding the corresponding algorithms. For example,     Chang and Yang \cite{w8} used a convex combination of all the previously generated primal iterative points   
    $$z_k=\frac{\psi-1}{\psi}x_{k-1}+\frac{1}{\psi}z_{k-1} $$
to replace $x_k$ in proximal term and presented   GRPDA. In GRPDA,  the extrapolation step \begin{equation*} \label{Exstep} \bar x_k=x_k+\theta(x_{k}-x_{k-1})
\end{equation*} in  PDA  is abandoned due to the introduction of the new point $z_k$.
    Though GRPDA usually shows good numerical results, it sometimes does not perform as well as the traditional PDA \cite{w15}. It may attribute to the abandon of  the extrapolation step $\bar x_k$ in \eqref{eq:PDA} of PDA.   
    
    It is universally acknowledged that the extrapolation step $\bar x_k $  
 plays an important role for convergence analysis and numerical performance of PDA. Therefore, it is interesting 
 to study some new  proximal terms  for  PDA without discarding the  extrapolation step $\bar x_k $.   
    Based on the above observation, we present a new primal-dual algorithm (NPDA), an  accelerated NPDA and    and   the NPDA with linesearch strategy. 
  The main contributions of the paper are summarized as below.
\begin{enumerate}
     \item By advantage of the previous iteration points, we introduce two auxiliary points $x_k^{ag}$ and $y_k^{ag}$, which are  a convex combination of all previous iterative points and are defined by 
     \begin{equation*} 
             x_{k}^{ag} =(1-a_k) x_{k-1}^{ag} + a_k x_{k-1},\;\;
            y_{k}^{ag} =(1-b_{k}) y_{k-1}^{ag} + b_{k} y_{k-1}.
        \end{equation*}     
The auxiliary points $x_k^{ag}$ and $y_k^{ag}$ are used    to generate two middle points:  
  \begin{equation*}
            x_{k}^{md} =(1-a_k) x_{k-1} + a_k x_{k}^{ag}, \;\;
             y_{k}^{md} =(1-b_{k}) y_{k-1} + b_{k} y_{k}^{ag}.
\end{equation*}
  
  We use $x_k^{md}$ and $y_k^{md}$ for the subproblems:
   \begin{align*}
            y_{k} &=\arg \max_{y \in Y} \{ \langle K \bar x_{k-1},y\rangle-f^*(y)-\frac{1}{2\sigma}\|y-y_k^{md}\|^2\},\\
            x_{k} &=\arg \min_{x \in X} \{g(x)+\langle Kx,y_k\rangle+\frac{1}{2\tau}\|x-x_k^{md}\|^2\},
        \end{align*}
  
  and present a new primal-dual algorithm for   problem \eqref{p}.
     The global convergence and $\mathcal{O}(1/N)$ ergodic rate of      NPDA are established. 

     
     \item When either $g$ or $f^*$ is strongly convex, we propose an adaptive strategy for parameters in NPDA and present an accelerated NPDA (ANPDA). We also establish   a  faster convergence rate $\mathcal{O}(1/N^2)$ for the primal variable $x$ or dual variable  $y$. 
     \item Since  the spectral norm of the linear transform $K$,  which is directly related to the conditions for stepsizes in NPDA, is   difficult to obtain in some applications, we introduce a linesearch strategy for NPDA and establish its convergence and  convergence rate.
     \item Numerical comparison with some state-of-the-art algorithms   on matrix game and LASSO problems demonstrate the promising performance of the proposed algorithms.
\end{enumerate}

\subsection{Organization.}
     The organization of the paper is as follows.   In Section \ref{section2}, we provide some notations, definitions and preliminary facts, which are crucial for the following theoretical analysis. Section \ref{section3} presents a  new primal-dual algorithm (NPDA), an accelerated NPDA  and NPDA with   linesearch.  The corresponding convergences  and convergence rates are also established. 
     Numerical experiments on  matrix game and LASSO problems are conducted  to evaluate the performance of the proposed algorithms in Section \ref{section4}. Finally, we  draw some concluding remarks in the last section.
\section{Preliminaries.}\label{section2}
   Some useful facts and lemmas are presented in the section, which will be often used  in the following text.

        \par  We first give the following two basic identities:       \begin{equation}
            2\langle a-b,a-c\rangle = \|a-b \|^2+\|a-c \|^2 -\| b-c \|^2,  \quad \forall a,b,c \in X.
        \end{equation}
        \begin{equation}
            \|\lambda x+(1-\lambda)y\|^2=\lambda\|x\|^2+(1-\lambda)\|y\|^2-\lambda(1-\lambda)\|x-y\|^2, \quad x,y \in X \;\; and \;\; \lambda \in R.
        \end{equation}
        
     \par If $(\hat{x},\hat{y})$ satisfies  
        \begin{align*}
            L(\hat{x},y)\leq L(\hat{x},\hat{y}) \leq L(x,\hat{y}), \quad \forall (x,y)\in X   \times Y,
        \end{align*}
    then   $(\hat{x},\hat{y})$ is  a saddle point of problem \eqref{p}. For any saddle point $(\hat{x},\hat{y})$,  it holds that $K \hat{x} \in \partial f^*(\hat{y})$ and $-K^T \hat{y} \in \partial g(\hat{x})$, where $\partial f^*$ and $\partial g$ denote the subgradients of  $f^*$ and $g$, respectively.  It also implies that   
        \begin{equation}
            P_{\hat{x},\hat{y}}(x)=g(x)-g(\hat{x})+\langle K^T\hat{y},x-\hat{x} \rangle \geq 0, \quad \forall x\in X,
        \end{equation}
        \begin{equation}
            D_{\hat{x},\hat{y}}(y)=f^*(y)- f^*(\hat{y})-\langle K\hat{x},y-\hat{y} \rangle \geq 0, \quad \forall y\in Y.
        \end{equation}
    The primal-dual gap function is defined by
        \begin{equation}\label{a5}
            \mathcal{G}_{\hat{x},\hat{y}}(x,y):= P_{\hat{x},\hat{y}}(x)+D_{\hat{x},\hat{y}}(y)\geq 0, \quad \forall (x,y)\in X   \times Y.
        \end{equation}
    Note that for any fixed saddle point $(\hat{x},\hat{y})$,  $P_{\hat{x},\hat{y}}(x)$, $D_{\hat{x},\hat{y}}(y)$ and $\mathcal{G}_{\hat{x},\hat{y}}(x,y)$ are convex.

    Let $h: X \rightarrow(-\infty,+\infty]$  be an extended  real-valued closed proper convex function and $\tau>0$. The proximal point mapping of $\tau h$ is defined  by
        \begin{align*}
            \text{Prox}_{\tau h}(x):=\arg \min_{y \in X} \{h(y)+\frac{1}{2\tau}\|y-x\|^2\}, \quad x \in X.
        \end{align*}
        
The following two lemmas give some properties of  the proximal point  map.    
    \begin{lemma}[\cite{w8}, Fact 2.1]
         For any  extended real-valued closed proper convex function $h: X\rightarrow(-\infty,+\infty]$, $\tau >0$  and $x\in X$, it holds that $z=\text{Prox}_{\tau h}(x)$ if and only if  $ \langle z-x,y-z   \rangle  \geq \tau ( h(z)-h(y) )$ for $\forall y\in X$, namely,         \begin{equation}\label{a1}
             z=\text{Prox}_{\tau h}(x) \;\Leftrightarrow \; \langle z-x,y-z\rangle \geq \tau ( h(z)-h(y) ), \quad \forall y\in X.
        \end{equation}
    \end{lemma}
    \begin{lemma}[\cite{w15}, Fact 2.1]\label{fact2.2}
        Let $h: X \rightarrow(-\infty,+\infty]$ be an extended real-valued closed proper and $\gamma$-strongly convex function with modulus $\gamma \geq 0$. Then for any $\tau > 0 $ and $x \in X$, it holds that $z = Prox_{\tau h}(x)$ if and only if $h(y) \geq h(z) + \frac{1}{\tau} \langle x - z, y - z\rangle +\frac{\gamma }{2}\| y - z\| ^2$ for all $y \in X$.
    \end{lemma}


\section{The proposed algorithms. }\label{section3}
  In the section, we  first present a new primal-dual algorithm for general convex-concave saddle point problems. We then develop  an accelerated version of the above proposed method   when  $f^*$ is strongly convex.  Finally, we incorporate a linesearch technique into the above  primal-dual algorithm.   The global convergence  and  ergodic convergence rate of the above three algorithms are also established.

\subsection{ New primal-dual algorithm (NPDA) for \eqref{p} when $g(x)$ and $f^*(y)$   are generally convex.   }
 We present a new primal-dual algorithm and establish its global convergence and ergodic convergence rate in the subsection.

Based on the motivation described in Section 1, the idea for the proposed algorithm   is  described as follows.

For   the affine equality constrained composite optimization problem:
\begin{align*}
    \min_{x \in X, w \in W}G(x)+F(w)  \\ s.t. \,\;\;\;\;   Bw-Kx=b,
  \end{align*}
where $X$ and $W$ are finite-dimensional Euclidean spaces, $G(\cdot)$ and $F(\cdot)$ are finitely valued, convex and lower
semi-continuous functions,   $K  \; \text{and} \; B  $   are bounded linear operators, the authors in \cite{w30}  proposed an accelerated   alternating direction method of multipliers (AADMM) to solve   the corresponding augmented Lagrangian formulation:
\begin{align*}
    \min_{x \in X, w \in W} \max_{y \in Y} G(x)+F(w)-\langle y, Bw-Kx-b \rangle +\frac{\rho}{2}\| Bw-Kx-b\|^2,
\end{align*}
where   $\rho >0 $ is a penalty parameter, $Y$ is finite-dimensional Euclidean spaces. Specifically, they  constructed auxiliary point $\bar x_{k}^{ag} =(1-a_{k-1}) \bar x_{k-1}^{ag} + a_{k-1} x_{k}$ and used a convex combination $\bar x_{k}^{md} =(1-a_{k}) \bar x_{k}^{ag} + a_{k} x_{k}$  for  the  subproblem about $x $:  
\begin{align*}
    x_{k+1}=&\arg \min_{x \in X}\{ \langle \nabla G(\bar x_k^{md}), x\rangle +\langle y_k,Kx\rangle+L_k+\frac{\eta_k}{2}\|x-x_k\|^2 \},  
\end{align*}
where $L_k=-\chi \theta_k \langle Bw_k-Kx_k-b,Kx\rangle+\frac{(1-\chi)\theta_k}{2}\|Bw_k-Kx-b\|^2$,    $ \theta_k>0,   \eta_k >0  $ and $   \chi=1 \; \text{or} \; 0$. It is noted that  $\bar w_k^{ag}$ and $\bar y_k^{ag}$ were also constructed by the way similarly to     $x_k^{ag}$ and   but not used   in the update steps of AADMM.  The introduction of $x_k^{md}$ in $\langle \nabla G(x_k^{md}), x\rangle$ can indeed accelerate   ADMM   \cite{w30}.

As described in Section 1, the proximal terms  $ \frac{1}{2\tau}\left\|x-x_k\right\|^2$ and $-\frac{1}{2\sigma}\left\|y-y_k\right\|^2$ are crucial to   PDA.  Motivated by \cite{w30}, we will construct two new points $x_k^{md}$ and $y_k^{md}$ (see Algorithm \ref{alg1}) and use them to  replace  $x_k $ and $y_k $ in $ \frac{1}{2\tau}\left\|x-x_k\right\|^2$ and $-\frac{1}{2\sigma}\left\|y-y_k\right\|^2$,  and present a new    primal-dual algorithm for solving problem \eqref{p} when $g$ and $f^*$ are both generally convex. The proposed algorithm   is described in detail as follows.
\begin{algorithm}[H] 
\caption{NPDA: The new primal-dual algorithm   for \eqref{p} when $g$ and $f $ are generally   convex    }  
\begin{algorithmic} 
    \State \textbf{Initialization:} Given two non-monotonic decreasing sequences $\{a_k\}$ and $\{b_k\}$, where $a_k, b_k \in (0,1)$, $\tau, \; \sigma>0, \; \theta =1$, $(x_0,y_1)\in X   \times Y$. Set $x_0^{ag}=x_0=x_{-1}$, $y_1^{ag}=y_1$, $\sqrt{\sigma \tau}L<1-a_k$ and $\sqrt{\sigma \tau}L<1-b_k$.
     \State \textbf{Main Iteration:}
     \State Step 1. Compute 
         \begin{equation}\label{Ax1al1}
             x_{k}^{ag} =(1-a_k) x_{k-1}^{ag} + a_k x_{k-1},
        \end{equation}
        \begin{equation}
            y_{k+1}^{ag} =(1-b_{k+1}) y_{k}^{ag} + b_{k+1} y_{k},
        \end{equation}
     \State   Step 2. Compute
        \begin{equation}\label{Ax2al1}
            x_{k}^{md} =(1-a_k) x_{k-1} + a_k x_{k}^{ag},
        \end{equation}
        \begin{equation}
             y_{k+1}^{md} =(1-b_{k+1}) y_k + b_{k+1} y_{k+1}^{ag},
        \end{equation}
    \State  Step 3. Compute       
         \begin{equation}
             \left \{ \begin{aligned}
             &x_{k}=\text{Prox}_{\tau g}(x_{k}^{md}-\tau K^Ty_{k}),\\
             &\bar{x}_{k}=x_{k}+\theta(x_{k}-x_{k-1}),\\
             &y_{k+1}=\text{Prox}_{\sigma f^*}(y_{k+1}^{md}+\sigma K\bar{x}_{k}).
        \end{aligned}\right.
        \end{equation}

    \State \textbf{End}
\end{algorithmic} \label{alg1}
\end{algorithm}

  \begin{remark}
      Since $a_k$ and $b_k$ in NPDA are     very small numbers, the conditions    $\sqrt{\sigma \tau}L<1-a_k$ and $\sqrt{\sigma \tau}L<1-b_k$  are not much stricter than the condition $\sqrt{\sigma \tau}L<1$ in PDA. 
 \end{remark}
\begin{remark}
  If $a_k  =  0$ and $b_k  =  0$ hold for $k\ge 1$ , then  NPAD reduces to the classical  primal-dual algorithm (PDA) \cite{w1}.
\end{remark}
\begin{remark}
    If $b_k=0$, $x_k^{md}=x_k^{ag}$ and $\theta=0$, then Algorithm 1 reduces to GRPDA \cite{w8}.
\end{remark}
\begin{remark} Though     $x_{k}^{md}$ in NPDA is motivated by \cite{w30}, there are important differences between NPDA and AADMM:    (i) In AADMM \cite{w30}, $\bar x_{k}^{md}$ is defined by $$\bar x_{k}^{md} =(1-a_{k}) \bar x_{k}^{ag} + a_{k} x_{k} ,\;\; \text{where} \;\;  \bar x_{k}^{ag} =(1-a_{k-1}) \bar x_{k-1}^{ag} + a_{k-1} x_{k},$$ while $x_{k}^{md}$ in   NPDA is given by $$x_{k}^{md} =a_k x_{k}^{ag} + (1-a_k) x_{k-1}   ,\;\;  \text{where} \;\;  x_{k}^{ag} =(1-a_k) x_{k-1}^{ag} + a_k x_{k-1}.$$ Obviously, the way for weights of $x_k^{md}$ and $x_k^{ag}$  are different from that for  $\bar x_k^{md}$ and  $\bar x_k^{ag}$. (ii)The use of $\bar x_k^{md}$ and $x_k^{md}$ is different. Specially,   $ \bar x_k^{md}$ in AADMM \cite{w30} is   used in  $ \langle \nabla G(\bar x_k^{md}), x\rangle$  rather than the proximal term,
   while $x_k^{md}$ in NPDA  is   used in the proximal  term $\frac{1}{2\tau}\|x-x_k^{md}\|^2$,  namely,    $$  x_{k}=  \arg \min_{x \in X}  \; \{g(x)+\langle Kx,y_k\rangle+\frac{1}{2\tau}\|x-x_k^{md}\|^2 \}. $$\\
    (iii)In NPDA,     $y_k^{md}$ and $y_k^{ag}$ are also used in the update of $y$,  while the corresponding term are not found in the update of $y$ or $w$ of AADMM \cite{w30}.
\end{remark}
    It is note that the use of   $x_k^{md}$  and $y_k^{md}$ in NPDA makes the convergence analysis quit complicated in contrast with that of PDA. We   establish the global convergence and ergodic convergence of NPDA  in the following two theorems.     
\begin{thm}[global convergence]\label{Th1}
 Suppose that  $\{(x_k,y_k)\}$ is the sequence generated by NPDA. Then, the sequence $\{(x_k,y_k)\}$ converges to a solution of problem (\ref{p}).
\end{thm} 

\begin{proof}
    Let $(\hat{x}, \hat{y})$ be any saddle point of problem (\ref{p}). It follow from (\ref{a1})  that 
       \begin{align}
            \langle x_{k+1}-x_{k+1}^{md} +\tau K^T y_{k+1},\hat{x}-x_{k+1}\rangle \; \geq \; \tau(g(x_{k+1})-g(\hat{x})),\label{p1}
       \end{align}
       \begin{align}
            \langle (y_{k+1}-y_{k+1}^{md})-\sigma K \bar{x}_k ,\hat{y}-y_{k+1}\rangle \; \geq \; \sigma (f^*(y_{k+1})-f^*(\hat{y})).\label{p2}
       \end{align}
   By summing   (\ref{p1}) and (\ref{p2}), we obtain that
        \begin{align}
            &\mathcal{G}_{\hat{x},\hat{y}}(x_{k+1},y_{k+1}) 
             +\langle K(x_{k+1}-\bar{x}_k),y_{k+1}-\hat{y} \rangle 
             +\frac{\| x_{k+1}-x_{k+1}^{md} \| ^2}{2\tau} \nonumber \\
            &\quad +\frac{\| x_{k+1}-\hat{x} \| ^2}{2\tau}
            +\frac{\| y_{k+1}-y_{k+1}^{md} \| ^2}{2\sigma}
             +\frac{\| y_{k+1}-\hat{y} \| ^2}{2\sigma} \nonumber \\
        &\leq \frac{\| x_{k+1}^{md}-\hat{x} \| ^2}{2\tau}+\frac{\| y_{k+1}^{md}-\hat{y} \| ^2}{2\sigma},\label{p3}
       \end{align}
   where $\mathcal{G}_{\hat{x},\hat{y}}(x_{k+1},y_{k+1})$ is given by \eqref{a5}.
  
    According to   \eqref{Ax1al1}  and  \eqref{Ax2al1}, we have
        \begin{align}
            \| x_{k+1}-\hat{x} \| ^2 =&\frac{1}{a_{k+2}}\| x_{k+2}^{ag}-\hat{x} \| ^2 -\frac{1-a_{k+2}}{a_{k+2}}\| x_{k+1}^{ag}-\hat{x} \| ^2+(1-a_{k+2})\| x_{k+1}-x_{k+1}^{ag} \|^2 ,
           \label{p4}
        \end{align}
        \begin{align}
            \| x_{k+1}-x_{k+1}^{md} \| ^2 =(1-a_{k+1}) \| x_{k+1}-x_k \| ^2 +a_{k+1} \| x_{k+1}^{ag}-x_{k+1} \| ^2 -a_{k+1} (1-a_{k+1})\| x_{k}-x_{k+1}^{ag} \| ^2, \label{p5}
        \end{align}
    and
        \begin{align}
            \| x_{k+1}^{md}-\hat{x} \| ^2 &= (1-a_{k+1}) \| x_k-\hat{x} \| ^2 + a_{k+1} \| x_{k+1}^{ag}-\hat{x} \| ^2 -a_{k+1}(1-a_{k+1})\| x_{k}-x_{k+1}^{ag} \| ^2 \nonumber \\
            &= \frac{1-a_{k+1}}{a_{k+1}} \| x_{k+1}^{ag}-\hat{x} \| ^2 -\frac{(1-a_{k+1})^2}{a_{k+1}}\| x_{k}^{ag}-\hat{x} \| ^2+(1-a_{k+1})^2\| x_{k}-x_{k}^{ag} \|^2 \nonumber \\
            &\quad + a_{k+1} \| x_{k+1}^{ag}-\hat{x} \| ^2 -a_{k+1}(1-a_{k+1})\| x_{k}-x_{k+1}^{ag} \| ^2.\label{p6}
        \end{align}
    We can also obtain  the  conclusions similar to (\ref{p4}), (\ref{p5}) and (\ref{p6}) for $y$.  
    
    By the definition of $\bar{x}_k$ and the fact that $\theta=1$, we can obtain that
        \begin{align}
            &\langle K(x_{k+1}-\bar{x}_k),y_{k+1}-\hat{y} \rangle \nonumber\\
            &=\langle K ( (x_{k+1}-x_k)-(x_k-x_{k-1})  ),y_{k+1}-\hat{y} \rangle \nonumber\\
            &=\langle K(x_{k+1}-x_k),y_{k+1}-\hat{y} \rangle-\langle K(x_{k}-x_{k-1}),y_{k}-\hat{y} \rangle \nonumber\\
            &- \langle K(x_{k}-x_{k-1}),y_{k+1}-y_{k}\rangle \nonumber\\
            &\geq  \langle K(x_{k+1}-x_k),y_{k+1}-\hat{y} \rangle- \langle K(x_{k}-x_{k-1}),y_{k}-\hat{y} \rangle \nonumber\\
            & - L\| x_{k}-x_{k-1} \| \| y_{k+1}-y_{k}\|.
        \end{align}
    \par For any $\alpha>0$, we have that (using $2ab \leq \alpha a^2+b^2/\alpha$ for any $a$, $b$)
        \begin{align}
            L\| x_{k}-x_{k-1} \| \| y_{k+1}-y_{k}\| \leq \frac{L\alpha \tau }{2\tau}\| x_{k}-x_{k-1} \| ^2 +\frac{L  \sigma}{2\alpha \sigma}\| y_{k+1}-y_{k}\|^2,\label{p8}
        \end{align}
    and if we set $\alpha =\sqrt{\sigma/\tau}$, then $L\alpha \tau = L\sigma/\alpha =\sqrt{\sigma \tau}L < 1$.
    \par By  (\ref{p3})-(\ref{p8}), we get
         \begin{align}
            &\mathcal{G}_{\hat{x},\hat{y}}(x_{k+1},y_{k+1})+\langle K(x_{k+1}-x_k),y_{k+1}-\hat{y} \rangle- \langle K(x_{k}-x_{k-1}),y_{k}-\hat{y} \rangle \nonumber\\
            &\quad +\frac{1}{2\tau}[(1-a_{k+1}) \| x_{k+1}-x_k \| ^2 +a_{k+1} \| x_{k+1}^{ag}-x_{k+1} \| ^2 -a_{k+1} (1-a_{k+1})\| x_{k}-x_{k+1}^{ag} \| ^2] \nonumber\\
            &\quad +\frac{1}{2\tau}[\frac{1}{a_{k+2}}\| x_{k+2}^{ag}-\hat{x} \| ^2 -\frac{1-a_{k+2}}{a_{k+2}}\| x_{k+1}^{ag}-\hat{x} \| ^2+(1-a_{k+2})\| x_{k+1}-x_{k+1}^{ag} \|^2 ] \nonumber\\
            &\quad +\frac{1}{2\sigma}[(1-b_{k+1}) \| y_{k+1}-y_k \| ^2 +b_{k+1} \| y_{k+1}^{ag}-y_{k+1} \| ^2 -b_{k+1} (1-b_{k+1})\| y_{k}-y_{k+1}^{ag} \| ^2] \nonumber\\
            &\quad +\frac{1}{2\sigma}[\frac{1}{b_{k+2}}\| y_{k+2}^{ag}-\hat{x} \| ^2 -\frac{1-b_{k+2}}{b_{k+2}}\| y_{k+1}^{ag}-\hat{x} \| ^2+(1-b_{k+2})\| y_{k+1}-y_{k+1}^{ag} \|^2 ] \nonumber\\
            &\leq
            \frac{1}{2\tau}[\frac{1-a_{k+1}}{a_{k+1}} \| x_{k+1}^{ag}-\hat{x} \| ^2 -\frac{(1-a_{k+1})^2}{a_{k+1}}\| x_{k}^{ag}-\hat{x} \| ^2+(1-a_{k+1})^2\| x_{k}-x_{k}^{ag} \|^2 \nonumber \\
            &\quad + a_{k+1} \| x_{k+1}^{ag}-\hat{x} \| ^2 -a_{k+1}(1-a_{k+1})\| x_{k}-x_{k+1}^{ag} \| ^2]\nonumber \\
             &\quad +\frac{1}{2\sigma}[\frac{1-b_{k+1}}{b_{k+1}} \| y_{k+1}^{ag}-\hat{x} \| ^2 -\frac{(1-b_{k+1})^2}{b_{k+1}}\| y_{k}^{ag}-\hat{x} \| ^2+(1-b_{k+1})^2\| y_{k}-y_{k}^{ag} \|^2 \nonumber \\
            &\quad +b_{k+1} \| y_{k+1}^{ag}-\hat{x} \| ^2 -b_{k+1}(1-b_{k+1})\| y_{k}-y_{k+1}^{ag} \| ^2]\nonumber \\
            &\quad +\frac{\sqrt{\sigma \tau}L}{2\tau}\|x_k-x_{k-1}\|^2+\frac{\sqrt{\sigma \tau}L}{2\sigma}\|y_{k+1}-y_{k}\|^2 .\label{p9}
        \end{align}
    Organizing inequality (\ref{p9}) yields
         \begin{align}
            &\mathcal{G}_{\hat{x},\hat{y}}(x_{k+1},y_{k+1})+\langle K(x_{k+1}-x_k),y_{k+1}-\hat{y} \rangle- \langle K(x_{k}-x_{k-1}),y_{k}-\hat{y} \rangle \nonumber\\
            &\quad +\frac{1}{2\tau}[(1-a_{k+1}) \| x_{k+1}-x_k \| ^2 +\left(1+a_{k+1}-a_{k+2}\right)\| x_{k+1}-x_{k+1}^{ag} \|^2
            +\frac{1}{a_{k+2}}\| x_{k+2}^{ag}-\hat{x} \| ^2 ] \nonumber\\
            &\quad +\frac{1}{2\sigma}[(1-b_{k+1}-\sqrt{\sigma \tau}L) \| y_{k+1}-y_k \| ^2 +\left(1+b_{k+1}-b_{k+2}\right)\| y_{k+1}-y_{k+1}^{ag} \|^2  +\frac{1}{b_{k+2}}\| y_{k+2}^{ag}-\hat{x} \| ^2]\nonumber\\
            &\leq
             \frac{1}{2\tau}[(\frac{1-a_{k+1}}{a_{k+1}}+\frac{1-a_{k+2}}{a_{k+2}}+a_{k+1}) \| x_{k+1}^{ag}-\hat{x} \| ^2 -\frac{(1-a_{k+1})^2}{a_{k+1}}\| x_{k}^{ag}-\hat{x} \| ^2+(1-a_{k+1})^2\| x_{k}-x_{k}^{ag} \|^2]\nonumber \\
             &\quad +\frac{1}{2\sigma}[\frac{1-b_{k+1}}{b_{k+1}}+ \frac{1-b_{k+2}}{b_{k+2}}+b_{k+1}) \| y_{k+1}^{ag}-\hat{x} \| ^2 -\frac{(1-b_{k+1})^2}{b_{k+1}}\| y_{k}^{ag}-\hat{x} \| ^2+(1-b_{k+1})^2\| y_{k}-y_{k}^{ag} \|^2 ]\nonumber \\
             &\quad +\frac{\sqrt{\sigma \tau}L}{2\tau}\|x_k-x_{k-1}\|^2 .\label{p10}
        \end{align}
    \par Summing (\ref{p10}) from $k=0$ to $N-1$ implies that   
         \begin{align}
            &\sum_{k=0}^{N-1}\mathcal{G}_{\hat{x},\hat{y}}(x_{k+1},y_{k+1})
            +\langle K(x_{N}-x_{N-1}),y_{N}-\hat{y} \rangle \nonumber\\
            &+\frac{1}{2\tau}[\sum_{k=0}^{N-2}(1-a_{k+1}-\sqrt{\sigma\tau}L) \| x_{k+1}-x_k \| ^2+(1-a_{N}) \| x_{N}-x_{N-1} \| ^2] \nonumber\\
            &+\frac{1}{2\tau}[\sum_{k=0}^{N-2}(a_{k+1}+a_{k+2}(1-a_{k+2})) \| x_{k+1}-x_{k+1}^{ag} \| ^2+(1+a_{N}-a_{N+1}) \| x_{N}-x_{N}^{ag} \| ^2] \nonumber\\
            &+\frac{1}{2\tau}[\sum_{k=0}^{N-3}(2-a_{k+2}-\frac{1}{a_{k+3}}+\frac{(1-a_{k+3})^2}{a_{k+3}}) \| x_{k+2}^{ag}-\hat{x} \| ^2\nonumber\\
            &+(2-a_{N}- \frac{1}{a_{N+1}}) \| x_{N}^{ag}-\hat{x} \| ^2+\frac{1}{a_{N+1}} \| x_{N+1}^{ag}-\hat{x} \| ^2] \nonumber\\
            &+\frac{1}{2\sigma}[\sum_{k=0}^{N-2}(1-b_{k+1}-\sqrt{\sigma\tau}L) \| y_{k+1}-y_k \| ^2+(1-b_{N}) \| y_{N}-y_{N-1} \| ^2] \nonumber\\
            &+\frac{1}{2\sigma}[\sum_{k=0}^{N-2}(b_{k+1}+b_{k+2}(1-b_{k+2})) \| y_{k+1}-y_{k+1}^{ag} \| ^2+(1+b_{N}-b_{N+1}) \| y_{N}-y_{N}^{ag} \| ^2] \nonumber\\
            &+\frac{1}{2\sigma}[\sum_{k=0}^{N-3}(2-b_{k+2}-\frac{1}{b_{k+3}}+\frac{(1-b_{k+3})^2}{b_{k+3}}) \| y_{k+2}^{ag}-\hat{y} \| ^2\nonumber\\
            &+(2-b_{N}- \frac{1}{b_{N+1}}) \| y_{N}^{ag}-\hat{y} \| ^2+\frac{1}{b_{N+1}} \| y_{N+1}^{ag}-\hat{y} \| ^2] \nonumber\\
            &\leq
            \frac{1}{2\tau}[(\frac{1}{a_1}+\frac{1}{a_2}+a_1-2-\frac{(1-a_2)^2}{a_2})\|x_1^{ag}-\hat{x}\|^2-\frac{(1-a_1)^2}{a_1}\|x_0^{ag}-\hat{x}\|^2] \nonumber\\
            &+\frac{1}{2\sigma}[(\frac{1}{b_1}+\frac{1}{b_2}+a_1-2-\frac{(1-b_2)^2}{b_2})\|y_1^{ag}-\hat{y}\|^2-\frac{(1-b_1)^2}{b_1}\|y_0^{ag}-\hat{y}\|^2]    \nonumber\\
            &=\frac{1}{2\tau}(2-a_2)\|x_0-\hat{x}\|^2
            +\frac{1}{2\sigma}(2-b_2)\|y_0-\hat{y}\|^2 \buildrel \Delta \over = C.  
         \label{p11}
        \end{align}
     By
        \begin{align*}
            2 \vert \langle K(x_{N}-x_{N-1}),y_N -\hat{y} \rangle  \vert
            \leq  \sqrt{\sigma \tau}L\| x_{N}-x_{N-1} \|^2/ \tau  +\sqrt{\sigma \tau}L
            \| y_{N}-\hat{y} \|^2/ \sigma,
        \end{align*}
     
        \begin{align}
           \frac{1}{a_{N+1}} \| x_{N+1}^{ag}-\hat{x} \| ^2=\frac{1-a_{N+1}}{a_{N+1}} \| x_{N}^{ag}-\hat{x} \| ^2+\| x_{N}-\hat{x} \| ^2-(1-a_{N+1})\| x_{N}-x_{N}^{ag} \| ^2 \label{X_N}
        \end{align}
        and
        \begin{align}
           \frac{1}{b_{N+1}} \| y_{N+1}^{ag}-\hat{y} \| ^2=\frac{1-b_{N+1}}{b_{N+1}} \| y_{N}^{ag}-\hat{x} \| ^2+\| y_{N}-\hat{y} \| ^2-(1-b_{N+1})\| y_{N}-y_{N}^{ag} \| ^2, \label{Y_N}
        \end{align}
       (\ref{p11}) can be rearranged as 
         \begin{align}
            &\frac{1}{2\tau}[\sum_{k=0}^{N-2}(1-a_{k+1}-\sqrt{\sigma\tau}L) \| x_{k+1}-x_k \| ^2+(1-a_{N}-\sqrt{\sigma\tau}L) \| x_{N}-x_{N-1} \| ^2] \nonumber\\
            &+\frac{1}{2\tau}[\sum_{k=0}^{N-2}(a_{k+1}+a_{k+2}(1-a_{k+2})) \| x_{k+1}-x_{k+1}^{ag} \| ^2+a_{N} \| x_{N}-x_{N}^{ag} \| ^2] \nonumber\\
            &+\frac{1}{2\tau}[\sum_{k=0}^{N-3}(2-a_{k+2}-\frac{1}{a_{k+3}}+\frac{(1-a_{k+3})^2}{a_{k+3}}) \| x_{k+2}^{ag}-\hat{x} \| ^2+(1-a_{N}) \| x_{N}^{ag}-\hat{x} \| ^2+ \| x_{N}-\hat{x} \| ^2]\nonumber\\
            &+\frac{1}{2\sigma}[\sum_{k=0}^{N-2}(1-b_{k+1}-\sqrt{\sigma\tau}L) \| y_{k+1}-y_k \| ^2+(1-b_{N}) \| y_{N}-y_{N-1} \| ^2] \nonumber\\
            &+\frac{1}{2\sigma}[\sum_{k=0}^{N-2}(b_{k+1}+b_{k+2}(1-b_{k+2})) \| y_{k+1}-y_{k+1}^{ag} \| ^2+b_{N} \| y_{N}-y_{N}^{ag} \| ^2] \nonumber\\
            &+\frac{1}{2\sigma}[\sum_{k=0}^{N-3}(2-b_{k+2}-\frac{1}{b_{k+3}}+\frac{(1-b_{k+3})^2}{b_{k+3}}) \| y_{k+2}^{ag}-\hat{y} \| ^2\nonumber\\
            &+(1-b_{N}) \| y_{N}^{ag}-\hat{y} \| ^2+ (1-\sqrt{\sigma\tau}L)\| y_{N}-\hat{y} \| ^2] +\sum_{k=0}^{N-1}\mathcal{G}_{\hat{x},\hat{y}}(x_{k+1},y_{k+1}) \nonumber\\
            &\leq C,
         \label{p12}
        \end{align}
        where $C$ is given by \eqref{p11}.
    Since the sequences $\{a_k\}$ and $\{b_k\}$ are non-monotonically  decreasing, it   deduces that the coefficients of each term in (\ref{p12}) are greater than or equal to 0.  Therefore, it follows that $\{x_k\} $ and $\{y_k\} $ are bounded, and 
         \begin{align*}
            \lim_{k \rightarrow +\infty} \|y_k -y_{k-1} \| ^2 =0,  \lim_{k \rightarrow +\infty} \|x_k -x_{k-1} \| ^2 =0,  \lim_{k \rightarrow +\infty}\| x_{k+1}^{ag}-x_{k+1} \|=0,
        \end{align*}
    which implies that 
        \begin{align*}
            x_{k+1}-x_{k+1}^{md}=[(1-a_{k+1})(x_{k+1}-x_{k})+
            a_{k+1}(x_{k+1}-x_{k+1}^{ag})] \rightarrow 0 \quad as \; k \rightarrow  +\infty,
        \end{align*}
    Similarly, 
        \begin{align*}
            y_{k+1}-y_{k+1}^{md} \rightarrow 0 \quad as \; k \rightarrow +\infty.
        \end{align*}
    Let $\{(x_{k_i},y_{k_i})\} $ be a subsequence that converges to a cluster point $(x^*,y^*)$, we deduce $x_{k_i}^{ag}$ and $x_{k_i}^{md}$ converges $x^*$, $y_{k_i}^{ag}$ and $y_{k_i}^{md}$ converges $y^*$. Taking the limit on two sides of the  following  two inequalities 
        \begin{align*}
            &\langle x_{{k_i}+1}-x_{{k_i}+1}^{md} +\tau K^T y_{{k_i}+1},x-x_{{k_i}+1}\rangle \; \geq \; \tau(g(x_{{k_i}+1})-g(x)), \quad \forall x\in X \\
            & \langle (y_{{k_i}+1}-y_{k_{i}+1}^{md})-\sigma K \bar{x}_{k_i} ,y-y_{{k_i}+1}\rangle \; \geq \; \sigma (f^*(y_{{k_i}+1})-f^*(y)), \quad \forall y\in Y
        \end{align*}
    we obtain that $(x^*,y^*)$ is a saddle point of (\ref{p}).
    \par Setting $(\hat{x},\hat{y})=(x^*,y^*)$ in (\ref{p12}), we have  $x_N \rightarrow x^*$, $y_N \rightarrow y^*$ as $N \rightarrow +\infty$.
\end{proof}

\begin{thm}[ergodic convergence]\label{Th2}
 Suppose that  $\{(x_k,y_k)\}$ is  generated by NPDA and $(\hat{x},\hat{y})$ is any saddle point of problem (\ref{p}). 
 Then, it holds that
 \begin{align*}
     \mathcal{G}_{\hat{x},\hat{y}}(X_N,Y_N) \leq
      \frac{C}{ N},
 \end{align*}
  where $X_N=\frac{1}{N}\sum_{k=1}^N x^k$, $Y_N=\frac{1}{N}\sum_{k=1}^N y^k$,  $C$ and  $\mathcal{G}_{\hat{x}, \hat{y}}(\cdot,\cdot)$ are given by   (\ref{p12}) and (\ref{a5}), respectively.
\end{thm} 

\begin{proof}
    Since  $(1-a_{k+1}-\sqrt{\sigma\tau}L)\geq0$, $(1-b_{k+1}-\sqrt{\sigma\tau}L)\geq0$ and $\{a_k\},\{b_k\}$ both are non-monotonically decreasing sequences,   by  (\ref{p12})  we have 
        \begin{align}
            \sum_{k=0}^{N-1}\mathcal{G}_{\hat{x},\hat{y}}(x_{k+1},y_{k+1})
            \leq C.
         \label{p13}
        \end{align} 
Together with the convexity of $g(\cdot)$ and $f^*(\cdot)$ and the definition of $X_N$ and $Y_N$, we can draw
        \begin{align*}
            \mathcal{G}_{\hat{x},\hat{y}}(X_{N},Y_{N})&=g(X_{N})-g(\hat{x})+\langle K^T\hat{y},X_{N}-\hat{x} \rangle+f^*(Y_{N})- f^*(\hat{y})-\langle K\hat{x},Y_{N}-\hat{y} \rangle \nonumber\\
            &=\frac{1}{N}\sum_{k=1}^{N}[g(X_{N})-g(\hat{x})+\langle K^T\hat{y},x_{k}-\hat{x} \rangle+f^*(Y_{N})- f^*(\hat{y})-\langle K\hat{x},y_{k}-\hat{y} \rangle] \nonumber\\
            &=\frac{1}{N}\sum_{k=0}^{N-1}[g(X_{N})-g(\hat{x})+\langle K^T\hat{y},x_{k+1}-\hat{x} \rangle+f^*(Y_{N})- f^*(\hat{y})-\langle K\hat{x},y_{k+1}-\hat{y} \rangle] \nonumber\\
            &\leq \frac{1}{N}\sum_{k=0}^{N-1}\mathcal{G}_{\hat{x},\hat{y}}(x_{k+1},y_{k+1})  \nonumber\\
           &\leq \frac{C}{N}.
        \end{align*}
\end{proof} 

\subsection{Accelerated NPDA for \eqref{p} when $f^*(y)$   is strongly convex.  }
     When either $g$ or $f^*$ in problem \eqref{p} is strongly convex, we present an accelerated NPDA (see   Algorithm \ref{alg2}) for solving problem \eqref{p}. 
     \par PDA often achieves  $\mathcal{O}(1/N)$ convergence rate \cite{w1}.   When either $g$ or $f^*$ in problem \eqref{p} is strongly convex, PDA can enjoy  faster $\mathcal{O}(1/N^2)$ convergence rate by choosing  stepsizes and $\theta_k$ adaptively. By adopting a similar technique, we can also accelerate NPDA when $f^*(y)$ in problem \eqref{p} is strongly convex. It is noted that      $f^*(y)$ is assumed to be strongly convex in the subsection, and $g$ is also assumed to be  strongly convex  due to the duality of primal-dual algorithms.

\begin{algorithm}[H] 
\caption{ANPDA: Accelerated NPDA   for problem \eqref{p} when $f^*$ is strongly convex.}  
\begin{algorithmic} 
    \State \textbf{Initialization:} Given two non-monotonic decreasing sequences $\{a_k\}$ and $\{b_k\}$, where $a_k, b_k \in (0,1)$,$\sigma_0>0$, $(x_0,y_1)\in X   \times Y$. Set $x_0^{ag}=x_0=x_{-1}$, $y_1^{ag}=y_1$, $\sqrt{\sigma \tau}L<1-a_k$ and $\sqrt{\sigma \tau}L<1-b_k$.
    \State \textbf{Main Iteration:}
    \State Step 1. Compute 
         \begin{equation}
             x_{k}^{ag} =(1-a_k) x_{k-1}^{ag} + a_k x_{k-1},
        \end{equation}
        \begin{equation}
            y_{k+1}^{ag} =(1-b_{k+1}) y_{k}^{ag} + b_{k+1} y_{k},
        \end{equation}
    \State   Step 2. Compute
        \begin{equation}
            x_{k}^{md} =(1-a_k) x_{k-1} + a_k x_{k}^{ag},
        \end{equation}
        \begin{equation}
             y_{k+1}^{md} =(1-b_{k+1}) y_k + b_{k+1} y_{k+1}^{ag},
        \end{equation}
    \State Step 3. Compute
        \begin{align*}
            x_{k}&=\text{Prox}_{\tau_{k-1} g}(x_{k}^{md}-\tau_{k-1} K^Ty_k),
        \end{align*}
    \State Step 4. Set 
        \begin{equation}
            \left \{
             \begin{aligned}\label{A}
               \theta_k&=1/{\sqrt{1+\gamma \sigma_{k-1}}},\\
                \sigma_{k}&=\sigma_{k-1} \theta_k,\\
                \tau_{k}&=\tau_{k-1} /\theta_{k},
            \end{aligned}\right.
      \end{equation}
    \State Step 5. Compute 
        \begin{align*}
             \bar{x}_{k}&=x_{k}+\theta_k(x_{k} - x_{k-1}),\\
            y_{k+1}&=\text{Prox}_{\sigma_k f^*}(y_{k+1}^{md}+\sigma_k K\bar{x}_{k}).\\           
        \end{align*}
    \State \textbf{End}       
\end{algorithmic}  \label{alg2}
\end{algorithm}  

The following two theorems present the main convergence results of ANPDA when $f^*$ is strong convex.  
\begin{thm}[global convergence]
  Suppose that $\{(x_k,y_k)\}$ is the sequence generated by ANPDA. Then, the sequence $\{(x_k,y_k)\}$ converges to a solution of problem (\ref{p}).
\end{thm} 

\begin{proof}
    It follows  from the strong convexity of $f^*$ and Lemma \ref{fact2.2} that
        \begin{align}
            \langle (y_{k+1}-y_{k+1}^{md})-\sigma_{k} K \bar{x}_k ,\hat{y}-y_{k+1}\rangle \; \geq \; \sigma_{k} (f^*(y_{k+1})-f^*(\hat{y})+\frac{\gamma}{2}\| y_{k+1}-\hat{y}\|^2).
        \end{align}
   Similarly to (\ref{p10}),  we obtain from the strong convexity of $f^*$ and  the above inequality that 
       \begin{align}
            &\mathcal{G}_{\hat{x},\hat{y}}(x_{k+1},y_{k+1})+\langle K(x_{k+1}-x_k),y_{k+1}-\hat{y} \rangle- \theta_k\langle K(x_{k}-x_{k-1}),y_{k}-\hat{y} \rangle \nonumber\\
            &\quad +\frac{1}{2\tau_{k}}[(1-a_{k+1}) \| x_{k+1}-x_k \| ^2 +\left(1+a_{k+1}-a_{k+2}\right)\| x_{k+1}-x_{k+1}^{ag} \|^2
            +\frac{1}{a_{k+2}}\| x_{k+2}^{ag}-\hat{x} \| ^2 ] \nonumber\\
            &\quad +\frac{1}{2\sigma_{k+1}}\frac{\sigma_{k+1}}{\sigma_{k}}(1+\gamma \sigma_{k})[(1-b_{k+1}-\sqrt{\sigma \tau}L) \| y_{k+1}-y_k \| ^2 \nonumber\\
            &\quad+\left(1+b_{k+1}-b_{k+2}\right)\| y_{k+1}-y_{k+1}^{ag} \|^2  +\frac{1}{b_{k+2}}\| y_{k+2}^{ag}-\hat{x} \| ^2]\nonumber\\
            &\leq
             \frac{1}{2\tau_{k}}[(\frac{1-a_{k+1}}{a_{k+1}}+\frac{1-a_{k+2}}{a_{k+2}}+a_{k+1}) \| x_{k+1}^{ag}-\hat{x} \| ^2 -\frac{(1-a_{k+1})^2}{a_{k+1}}\| x_{k}^{ag}-\hat{x} \| ^2+(1-a_{k+1})^2\| x_{k}-x_{k}^{ag} \|^2]\nonumber \\
            &\quad +\frac{1}{2\sigma_{k}}[\frac{1-b_{k+1}}{b_{k+1}}+ \frac{1-b_{k+2}}{b_{k+2}}+b_{k+1}) \| y_{k+1}^{ag}-\hat{x} \| ^2 -\frac{(1-b_{k+1})^2}{b_{k+1}}\| y_{k}^{ag}-\hat{x} \| ^2+(1-b_{k+1})^2\| y_{k}-y_{k}^{ag} \|^2 ]\nonumber \\
            &\quad +\frac{\sqrt{\sigma_{k} \tau_{k}}L}{2\tau_{k}}\|x_k-x_{k-1}\|^2 \label{pa1}.  
       \end{align}
   By the definitions of $\theta_{k}$, $\tau_{k}$ and $\sigma_{k}$, we  get
       \begin{align*}
           (1+\gamma \sigma_{k})\frac{\sigma_{k+1}}{\sigma_{k}}=\frac{\tau_{k+1}}{\tau_{k}}=\frac{1}{\theta_{k}}>1 
       \end{align*}
    and $$\sigma_{k} \tau_{k}=\sigma_{0} \tau_{0}.$$
     Multiplying  (\ref{pa1})  by $\frac{1}{\sigma_k }$ yields
       \begin{align}
            &\frac{1}{\sigma_k}\mathcal{G}_{\hat{x},\hat{y}}(x_{k+1},y_{k+1})+\frac{1}{\sigma_k}\langle K(x_{k+1}-x_k),y_{k+1}-\hat{y} \rangle- \frac{1}{\sigma_{k-1}}\langle K(x_{k}-x_{k-1}),y_{k}-\hat{y} \rangle \nonumber\\
            &\quad +\frac{1}{2\tau_{0}\sigma_0}[(1-a_{k+1}) \| x_{k+1}-x_k \| ^2 +\left(1+a_{k+1}-a_{k+2}\right)\| x_{k+1}-x_{k+1}^{ag} \|^2
            +\frac{1}{a_{k+2}}\| x_{k+2}^{ag}-\hat{x} \| ^2 ] \nonumber\\
            &\quad +\frac{1}{2\sigma_{k+1}\sigma_{k+1}}[(1-b_{k+1}-\sqrt{\sigma \tau}L) \| y_{k+1}-y_k \| ^2+\left(1+b_{k+1}-b_{k+2}\right)\| y_{k+1}-y_{k+1}^{ag} \|^2  +\frac{1}{b_{k+2}}\| y_{k+2}^{ag}-\hat{x} \| ^2]\nonumber\\
            &\leq
             \frac{1}{2\tau_{0}\sigma_0}[(\frac{1-a_{k+1}}{a_{k+1}}+\frac{1-a_{k+2}}{a_{k+2}}+a_{k+1}) \| x_{k+1}^{ag}-\hat{x} \| ^2 -\frac{(1-a_{k+1})^2}{a_{k+1}}\| x_{k}^{ag}-\hat{x} \| ^2+(1-a_{k+1})^2\| x_{k}-x_{k}^{ag} \|^2]\nonumber \\
            &\quad +\frac{1}{2\sigma_{k}\sigma_k}[(\frac{1-b_{k+1}}{b_{k+1}}+ \frac{1-b_{k+2}}{b_{k+2}}+b_{k+1}) \| y_{k+1}^{ag}-\hat{x} \| ^2 -\frac{(1-b_{k+1})^2}{b_{k+1}}\| y_{k}^{ag}-\hat{x} \| ^2+(1-b_{k+1})^2\| y_{k}-y_{k}^{ag} \|^2 ]\nonumber \\
            &\quad +\frac{\sqrt{\sigma_{0} \tau_{0}}L}{2\tau_{0}\sigma_0}\|x_k-x_{k-1}\|^2 \label{pa2}.  
       \end{align} 
    Summing (\ref{pa2}) from $k=0$ to $N-1$ and using the inequality
        \begin{align*}
            2 \vert \langle K(x_{N}-x_{N-1}),y_N -\hat{y} \rangle  \vert
            \leq  \sqrt{\sigma_{N-1} \tau_{N-1}}L\| x_{N}-x_{N-1} \|^2/ \tau_{N-1}  +\sqrt{\sigma_{N-1} \tau_{N-1}}L \| y_{N}-\hat{y} \|^2/ \sigma_{N-1},
        \end{align*}
    we obtain  that 
     
        \begin{align} \label{i42}
            &\sum_{k=0}^{N-1}\frac{1}{\sigma_k}\mathcal{G}_{\hat{x},\hat{y}}(x_{k+1},y_{k+1})+\frac{\tau_0 \sigma_0}{2\sigma_{N-1}^2} \| y_{N}-\hat{y} \| ^2 \nonumber\\
            &+\frac{1}{2\tau_0 \sigma_0}[\sum_{k=0}^{N-2}(1-a_{k+1}-\sqrt{\tau_0 \sigma_0}L) \| x_{k+1}-x_k \| ^2+(1-a_{N+1}-\sqrt{\tau_0 \sigma_0}L) \| x_{N}-x_{N-1} \| ^2] \nonumber\\
            &+\frac{1}{2\tau_0 \sigma_0}[\sum_{k=0}^{N-2}(a_{k+1}+a_{k+2}(1-a_{k+2})) \| x_{k+1}-x_{k+1}^{ag} \| ^2+(1+a_{N}-a_{N+1}) \| x_{N}-x_{N}^{ag} \| ^2] \nonumber\\
            &+\frac{1}{2\tau_0 \sigma_0}[\sum_{k=0}^{N-3}(2-a_{k+2}-\frac{1}{a_{k+3}}+\frac{(1-a_{k+3})^2}{a_{k+3}}) \| x_{k+2}^{ag}-\hat{x} \| ^2\nonumber\\
            &+(2-a_{N}- \frac{1}{a_{N+1}}) \| x_{N}^{ag}-\hat{x} \| ^2+\frac{1}{a_{N+1}} \| x_{N+1}^{ag}-\hat{x} \| ^2] \nonumber\\
            &+\sum_{k=0}^{N-2}\frac{1}{2\sigma_{k+1}^2}(1-b_{k+1}-\sqrt{\tau_0 \sigma_0}L) \| y_{k+1}-y_k \| ^2+\frac{1}{2\sigma_{N}^2}(1-b_{N}) \| y_{N}-y_{N-1} \| ^2 \nonumber\\
            &+\sum_{k=0}^{N-2}\frac{1}{2\sigma_{k+1}^2}(b_{k+1}+b_{k+2}(1-b_{k+2})) \| y_{k+1}-y_{k+1}^{ag} \| ^2+\frac{1}{2\sigma_{N}^2}(1+b_{N}-b_{N+1}) \| y_{N}-y_{N}^{ag} \| ^2 \nonumber\\
            &+\sum_{k=0}^{N-3}\frac{1}{2\sigma_{k+1}^2}(2-b_{k+2}-\frac{1}{b_{k+3}}+\frac{(1-b_{k+3})^2}{b_{k+3}}) \| y_{k+2}^{ag}-\hat{y} \| ^2\nonumber\\
            &+\frac{1}{2\sigma_{N-1}^2}(2-b_{N}- \frac{1}{b_{N+1}}) \| y_{N}^{ag}-\hat{y} \| ^2+\frac{1}{2\sigma_{N}^2}\frac{1}{b_{N+1}} \| y_{N+1}^{ag}-\hat{y} \| ^2 \nonumber\\
            &\leq
            \frac{1}{2\tau_0 \sigma_0}[(\frac{1}{a_1}+\frac{1}{a_2}+a_1-2-\frac{(1-a_2)^2}{a_2})\|x_1^{ag}-\hat{x}\|^2-\frac{(1-a_1)^2}{a_1}\|x_0^{ag}-\hat{x}\|^2] \nonumber\\
            &+\frac{1}{2\sigma_0^2}[(\frac{1}{b_1}+\frac{1}{b_2}+a_1-2-\frac{(1-b_2)^2}{b_2})\|y_1^{ag}-\hat{y}\|^2-\frac{(1-b_1)^2}{b_1}\|y_0^{ag}-\hat{y}\|^2]\nonumber\\
            &=\frac{1}{2\tau_0 \sigma_0}(2-a_2)\|x_0-\hat{x}\|^2+\frac{1}{2\sigma_0^2}(2-b_2)\|y_0-\hat{y}\|^2 \buildrel \Delta \over = C_1. 
        \end{align}
 By   using $\frac{1}{\sigma_{N-1}} \leq \frac{1}{\sigma_{N}}$, (\ref{X_N}) and (\ref{Y_N}), we get    
        \begin{align}
            &\frac{1}{2\tau_0 \sigma_0}[\sum_{k=0}^{N-2}(1-a_{k+1}-\sqrt{\tau_0 \sigma_0}L) \| x_{k+1}-x_k \| ^2+(1-a_{N+1}-\sqrt{\tau_0 \sigma_0}L) \| x_{N}-x_{N-1} \| ^2] \nonumber\\
            &+\frac{1}{2\tau_0 \sigma_0}[\sum_{k=0}^{N-2}(a_{k+1}+a_{k+2}(1-a_{k+2})) \| x_{k+1}-x_{k+1}^{ag} \| ^2+a_{N} \| x_{N}-x_{N}^{ag} \| ^2] \nonumber\\
            &+\frac{1}{2\tau_0 \sigma_0}[\sum_{k=0}^{N-3}(2-a_{k+2}-\frac{1}{a_{k+3}}+\frac{(1-a_{k+3})^2}{a_{k+3}}) \| x_{k+2}^{ag}-\hat{x} \| ^2+(1-a_{N}) \| x_{N}^{ag}-\hat{x} \| ^2+ \| x_{N}-\hat{x} \| ^2]\nonumber\\
            &+\sum_{k=0}^{N-2}\frac{1}{2\sigma_{k+1}^2}(1-b_{k+1}-\sqrt{\tau_0 \sigma_0}L) \| y_{k+1}-y_k \| ^2+\frac{1}{2\sigma_{N}^2}(1-b_{N}) \| y_{N}-y_{N-1} \| ^2 \nonumber\\
            &+\sum_{k=0}^{N-2}\frac{1}{2\sigma_{k+1}^2}(b_{k+1}+b_{k+2}(1-b_{k+2})) \| y_{k+1}-y_{k+1}^{ag} \| ^2+\frac{1}{2\sigma_{N}^2}b_{N} \| y_{N}-y_{N}^{ag} \| ^2] \nonumber\\
            &+\sum_{k=0}^{N-3}\frac{1}{2\sigma_{k+1}^2}(2-b_{k+2}-\frac{1}{b_{k+3}}+\frac{(1-b_{k+3})^2}{b_{k+3}}) \| y_{k+2}^{ag}-\hat{y} \| ^2\nonumber\\
            &+\frac{1}{2\sigma_{N-1}^2}(1-b_{N}) \| y_{N}^{ag}-\hat{y} \| ^2+\frac{1}{2\sigma_{N-1}^2} (1-\sqrt{\tau_0 \sigma_0}L)\| y_{N}-\hat{y} \| ^2 +\sum_{k=0}^{N-1}\frac{1}{\sigma_k}\mathcal{G}_{\hat{x},\hat{y}}(x_{k+1},y_{k+1})\nonumber\\
            &\leq
            C_1,\label{pa3}  
        \end{align}
    where  $C_1$ is given by \eqref{i42}. 
    
    Similarly to  Theorem \ref{Th1}, we can obtain the global convergence of ANPDA.                
\end{proof}

\begin{thm}[convergence rate]
    Suppose that $\{(x_k,y_k)\}$ is  generated by ANPDA and $(\hat{x},\hat{y})$ is any saddle point of problem (\ref{p}). Then, it hold that
        \begin{align*}
            \|y_N-\hat{y}\|=\mathcal{O}(\frac{1}{N^2}),
        \end{align*}
    and
        \begin{align*}
           \mathcal{G}_{\hat{x},\hat{y}}(X_N,Y_N) \leq
            \frac{C_1}{N}, 
        \end{align*}
    where $X_N=\frac{1}{N}\sum_{k=1}^N x^k$, $Y_N=\frac{1}{N}\sum_{k=1}^N y^k$, and  $C_1$ is given by \eqref{i42}. 
\end{thm} 

\begin{proof}
    It follows from  (\ref{pa3}) that  
    $(1-\sqrt{\tau_0 \sigma_0}L)\| y_N-\hat{y} \|^2
    \leq 2\sigma_{N-1}^2 C_1$, which  together with $\gamma \sigma_N \sim N^{-1}$ ( Corollary 1 in \cite{w1}) implies that $\|y_N-\hat{y}\|=\mathcal{O}(\frac{1}{N^2})$.
    \par   Let $X_N=\frac{1}{N}\sum_{k=1}^N x_k$, $Y_N=\frac{1}{N}\sum_{k=1}^N y_k$, and $1/\sigma_k \geq 1/\sigma_0$ for $  k\geq 0$ in $\mathcal{G}_{\hat{x},\hat{y}}(x,y)$. Similarly to Theorem \ref{Th2}, we can draw that  $ \mathcal{G}_{\hat{x},\hat{y}}(X_N,Y_N) 
    \leq
    \frac{1}{ N}\sigma_0 C_1$.
\end{proof}

\subsection{NPDA with linesearch.}  
  Since the conditions $\sqrt{\sigma \tau}L<1-a_k$ and $\sqrt{\sigma \tau}L<1-b_k$ for stepsizes in NPDA is related directly to the spectral norm $L$ of the linear transform,   which is difficult to obtain in some applications, 
  we   propose a linsearch strategy for  NPDA    and establish its global convergence and convergence rate in the subsection. 
\begin{algorithm}[H]
\caption{NPDAL: NPDA with linesearch  }
\begin{algorithmic} 
    \State \textbf{Initialization:} Given two non-monotonic decreasing sequences $\{a_k\}$ and $\{b_k\}$, where $a_k, b_k \in (0,1)$, $(x_0,y_1)\in X   \times Y$, $\tau_0>0$, $\theta>0$, $\mu \in (0,1)$, $\delta \in (0,1)$ with $\delta+\theta^2a_k+a_k<1$ and $\delta+b_k<1$, $\beta>0$. Set $x_0^{ag}=x_0$, $y_1^{ag}=y_1$, $\theta_0=1$.
    \State \textbf{Main Iteration:}
         \State Step 1. Compute 
         \begin{equation}
             x_{k}^{ag} =(1-a_k) x_{k-1}^{ag} + a_k x_{k-1},
        \end{equation}
        \begin{equation}
            y_{k+1}^{ag} =(1-b_{k+1}) y_{k}^{ag} + b_{k+1} y_{k},
        \end{equation}
     \State   Step 2. Compute
        \begin{equation}
            x_{k}^{md} =(1-a_k) x_{k-1} + a_k x_{k}^{ag},
        \end{equation}
        \begin{equation}
             y_{k+1}^{md} =(1-b_{k+1}) y_k + b_{k+1} y_{k+1}^{ag},
        \end{equation}
    \State Step 3. Compute
         \begin{align}
            x_{k}&=\text{Prox}_{\tau_{k-1} g}(x_{k}^{md}-\tau_{k-1} K^Ty_k),\label{algl1}
        \end{align}
    \State Step 4. Set $\tau_{k} =\tau_{k-1} \sqrt{1+\theta_{k-1}}$ and run  the following \textbf{ linesearch:}
    \State \hspace{3.2em} 4.a. Compute
        \begin{align}
            \theta_{k}&=\frac{\tau_{k}}{\tau_{k-1}},\\
            \bar{x}_{k}&=x_{k}+\theta_k(x_{k} - x_{k-1}),\\
            y_{k+1}&=\text{Prox}_{\beta\tau_k f^*}(y_{k+1}^{md}+\beta\tau_k K\bar{x}_k),
        \end{align}
    \State  \hspace{3.2em} 4.b. where set $\tau_{k} :=\mu\tau_{k}  $ until  the following condition are met
        \begin{align*}
            \sqrt{\beta}\tau_k||K^Ty_{k+1}-K^Ty_{k}||\leq\delta ||y_{k+1}-y_{k}||.
        \end{align*}       
    \State  \hspace{3.2em} \textbf{End of linesearch}
    \State \textbf{End} 
\end{algorithmic}  \label{alg3}
\end{algorithm}  

\begin{remark}
   In NPDAL, $\theta$ is set to  an upper bound for $\theta_k$, for example, $\theta=\frac{\sqrt{5}+1}{2}$ (see Lemma \ref{lemma3.1} and \cite{w2}).
\end{remark}
\begin{remark}\label{L2}
      It is require to calculate $\text{Prox}_{\beta\tau_k f^* }(\cdot)$ and $K^T y_{k+1}$ in the linesearch of NPDAL. As  pointed out in \cite{w2},   this procedure becomes extremely simple when $\text{Prox}_{\sigma f^*}$ is linear operator or affine operator. Here are some simple examples:
\end{remark}
\begin{enumerate}[label=(\alph*)]
     \item $\text{Prox}_{\lambda f^*} (u) = u - \lambda b $ when $f^* (y) = \langle b, y\rangle$ for some $b \in Y $;
     \item $\text{Prox}_{\lambda f^*} (u) = \frac{1}{1+\lambda}(u + \lambda b)$ when $f^* (y) = \frac{1}{2}\| y - b\| ^2$ for some $b \in Y $;
     \item $\text{Prox}_{\lambda f^*} (u)  = u +\frac{b - \langle a,u\rangle}{\| a\| ^2}a$ when $f^* (y)$ is the indicator function of $H = \{ u :\langle a, u\rangle = b\}$ for some  $a \in Y $ and  $b \in R $.
\end{enumerate}
For the linesearch in NPDAL, we can obtain  the following lemma.
\begin{lemma}  \label{lemma3.1}
     (i) The linesearch in NPDAL always terminates.\\
    (ii) There exists $\tau> 0$ such that $\tau_k > \tau$ for all $k\geq 0$.\\
    (iii) There exists $\theta > 0$ such that $\theta_k \leq \theta$ for all $k \geq 0$.
\end{lemma}
\begin{proof}
    The detailed proof is given in Lemma 3.3 of  \cite{w2}.
\end{proof}
The following two theorems  present    the global convergence and convergence rate of NPDAL.
\begin{thm}[global convergence]
    Suppose that  $\{(x_k,y_k)\}$ is the sequence generated   by NPDAL. Then,     the sequence $\{(x_k,y_k)\}$ converges to a solution of problem (1).
\end{thm} 

\begin{proof}
    Let $(\hat{x}, \hat{y})$ be any saddle point of problem (\ref{p}). It follows from  (\ref{a1})  that
        \begin{align}
            \langle x_{k+1}-x_{k+1}^{md} +\tau_k K^T y_{k+1},\hat{x}-x_{k+1}\rangle \; \geq \; \tau_k \left(g(x_{k+1})-g(\hat{x})\right),
            \label{l1}
        \end{align}
        \begin{align}
            \langle \frac{1}{\beta}(y_{k+1}-y_{k+1}^{md})-\tau_k K \bar{x}_k ,\hat{y}-y_{k+1}\rangle \; \geq \; \tau_k (f^*(y_{k+1})-f^*(\hat{y})).
            \label{l2}
        \end{align}
   By Step 3 of NPDAL and (\ref{a1}), we obtain  that  
        \begin{align*}
            \langle x_{k}-x_{k}^{md} +\tau_{k-1} K^T y_{k},x-x_{k}\rangle \; \geq \; \tau_{k-1}(g(x_{k})-g(x)), \;\; \forall x \in X.
        \end{align*}
    Setting $x=x_{k+1}$ and $x=x_{k-1}$ in the above inequality, we get
        \begin{align}
            \langle x_{k}-x_{k}^{md} +\tau_{k-1} K^T y_{k},x_{k+1}-x_{k}\rangle \; \geq \; \tau_{k-1}(g(x_{k})-g(x_{k+1})),
            \label{l3}
        \end{align}
        \begin{align}
            \langle x_{k}-x_{k}^{md} +\tau_{k-1} K^T y_{k},x_{k-1}-x_{k}\rangle \; \geq \; \tau_{k-1}(g(x_{k})-g(x_{k-1})).
            \label{l4}
        \end{align}
    Adding (\ref{l3}), multiplied by $\theta_k$ and (\ref{l4}), multiplied by $\theta_k^2$, and use $\bar{x}_{k}=x_{k}+\theta_k(x_{k} - x_{k-1})$, we obtain
        \begin{align}
             \langle \theta_k( x_{k}-x_{k}^{md}) +\tau_{k} K^T y_{k},x_{k+1}-\bar{x}_{k}\rangle \;
            \geq \tau_{k}((1+\theta_k)g(x_k)-g(x_{k+1})-\theta_k g(x_{k-1})).\label{l5}
        \end{align}
    Summing (\ref{l1}), (\ref{l2}) and (\ref{l5}) yields that
        \begin{align}
            &\langle x_{k+1} - x_{k+1}^{md}, \hat{x} - x_{k+1} \rangle + \frac{1}{\beta} \langle y_{k+1} - y_{k+1}^{md}, \hat{y} - y_{k+1} \rangle  + \langle \theta_k( x_{k}-x_{k}^{md}), x_{k+1} - \bar{x}_{k} \rangle \nonumber\\
            &\quad + \tau_{k} \langle K^T y_{k}, x_{k+1} - \bar{x}_{k} \rangle- \tau_{k} \langle K \bar{x}_k, \hat{y} - y_{k+1} \rangle + \tau_{k} \langle K^T y_{k+1}, \hat{x} - x_{k+1} \rangle \nonumber\\
            &\geq \tau_k \left( f^*(y_{k+1}) - f^*(\hat{y}) + (1 + \theta_k) g(x_k) - g(x_{k-1}) - \theta_k g(\hat{x}) \right).
        \end{align} 
    Applying the definition of $D_{\hat{x},\hat{y}}(y_{k+1})$, $P_{\hat{x},\hat{y}}(x_{k+1})$, and adding $\tau_k\langle K^T \hat{y},\bar{x}_k-\hat{x} \rangle-\langle K\hat{x},y_{k+1}-\hat{y}\rangle$ to both sides of the above inequality  simultaneously, we get
        \begin{align}
             &\langle x_{k+1} - x_{k+1}^{md}, \hat{x} - x_{k+1} \rangle + \frac{1}{\beta} \langle y_{k+1} - y_{k+1}^{md}, \hat{y} - y_{k+1} \rangle  + \langle \theta_k( x_{k}-x_{k}^{md}), x_{k+1} - \bar{x}_{k} \rangle \nonumber\\
             &+\tau_k \langle K^T  y_{k+1}-K^T y_{k},\bar{x}_{k}-x_{k+1}\rangle  \nonumber\\
             &\geq \tau_k ( (1+\theta_k)P_{\hat{x},\hat{y}}(x_{k})-\theta_k P_{\hat{x},\hat{y}}(x_{k-1})+D_{\hat{x},\hat{y}}(y_{k+1}) ) \buildrel \Delta \over =  \omega_k .\label{l6}
        \end{align}
    Applying the cosine rule the left side of  (\ref{l6})  and using 
        \begin{align*}
            2\tau_k\langle K^T  y_{k+1}-K^T y_{k},\bar{x}_{k}-x_{k+1}\rangle 
             &\leq  \frac{2\delta}{\sqrt{\beta}}\|\bar{x}_{k}-x_{k+1}\| \| y_{k+1}-y_{k}\| \\
             &\leq  \delta\|\bar{x}_{k}-x_{k+1}\|^2+\frac{\delta}{\beta}\| y_{k+1}- y_{k}\|^2 
        \end{align*}
    and
        \begin{align}
             &{ 2 \langle \theta_k( x_{k}-x_{k}^{md}), x_{k+1} - \bar{x}_{k} \rangle } \nonumber\\
             &=2 (1-a_k)\theta_k\langle  x_{k}-x_{k-1},  x_{k+1} - \bar{x}_{k}  \rangle+ 2\theta_k a_k\langle  x_{k}-x_{k}^{ag}, x_{k+1} - \bar{x}_{k} \rangle \nonumber\\
              &=2(1-a_k)\langle  \bar{x}_{k}-x_{k},  x_{k+1} - \bar{x}_{k}  \rangle +2\theta_k a_k\langle  x_{k}-x_{k}^{ag}, x_{k+1} - \bar{x}_{k} \rangle \nonumber\\
              &\leq (1-a_k)\|x_{k+1}-x_{k}\|^2-(1-a_k)\|\bar{x}_{k}-x_{k}\|^2-(1-a_k)\|x_{k+1}-\bar{x}_{k}\|^2\nonumber\\
              &+ a_k\| x_{k}-x_{k}^{ag}\|^2+\theta_k^2 a_k\|x_{k+1} - \bar{x}_{k}\|^2,
        \end{align}
    we have
        \begin{align}
           & 2\omega_k +\| x_{k+1}-x_{k+1}^{md} \| ^2 +\| x_{k+1}-\hat{x} \| ^2  
           +\frac{1}{\beta}\| y_{k+1}-y_{k+1}^{md} \| ^2 +\frac{1}{\beta}\| y_{k+1}-\hat{y} \| ^2 \nonumber \\
           &\leq \| x_{k+1}^{md}-\hat{x} \| ^2+\frac{1}{\beta}\| y_{k+1}^{md}-\hat{y} \| ^2 -(1-a_k-\theta_k^2 a_k-\delta)\|x_{k+1}-\bar{x}_{k}\|^2+\frac{\delta}{\beta}\| y_{k+1}- y_{k}\|^2 \nonumber\\
           &+(1-a_k)\|x_{k+1}-x_{k}\|^2-(1-a_k)\|\bar{x}_{k}-x_{k}\|^2 + a_k\| x_{k}-x_{k}^{ag}\|^2.\label{l7}
        \end{align}
    Substituting   (\ref{p4}),  (\ref{p6}) and the corresponding results for $y$   into (\ref{l7})  yields
        \begin{align}
            &2\omega_k 
             +[a_{k+1} \| x_{k+1}^{ag}-x_{k+1} \| ^2 -a_{k+1} (1-a_{k+1})\| x_{k}-x_{k+1}^{ag} \| ^2] \nonumber\\
            &\quad +[\frac{1}{a_{k+2}}\| x_{k+2}^{ag}-\hat{x} \| ^2 -\frac{1-a_{k+2}}{a_{k+2}}\| x_{k+1}^{ag}-\hat{x} \| ^2+(1-a_{k+2})\| x_{k+1}-x_{k+1}^{ag} \|^2 ] \nonumber\\
            &\quad +\frac{1}{\beta}[(1-b_{k+1}-\delta) \| y_{k+1}-y_k \| ^2 +b_{k+1} \| y_{k+1}^{ag}-y_{k+1} \| ^2 -b_{k+1} (1-b_{k+1})\| y_{k}-y_{k+1}^{ag} \| ^2] \nonumber\\
            &\quad +\frac{1}{\beta}[\frac{1}{b_{k+2}}\| y_{k+2}^{ag}-\hat{x} \| ^2 -\frac{1-b_{k+2}}{b_{k+2}}\| y_{k+1}^{ag}-\hat{x} \| ^2+(1-b_{k+2})\| y_{k+1}-y_{k+1}^{ag} \|^2 ] \nonumber\\
            &\leq
            [\frac{1-a_{k+1}}{a_{k+1}} \| x_{k+1}^{ag}-\hat{x} \| ^2 -\frac{(1-a_{k+1})^2}{a_{k+1}}\| x_{k}^{ag}-\hat{x} \| ^2+((1-a_{k+1})^2+a_k)\| x_{k}-x_{k}^{ag} \|^2 \nonumber \\
            &\quad + a_{k+1} \| x_{k+1}^{ag}-\hat{x} \| ^2 -a_{k+1}(1-a_{k+1})\| x_{k}-x_{k+1}^{ag} \| ^2]\nonumber \\
             &\quad +\frac{1}{\beta}[\frac{1-b_{k+1}}{b_{k+1}} \| y_{k+1}^{ag}-\hat{x} \| ^2 -\frac{(1-b_{k+1})^2}{b_{k+1}}\| y_{k}^{ag}-\hat{x} \| ^2+(1-b_{k+1})^2\| y_{k}-y_{k}^{ag} \|^2 \nonumber \\
            &\quad +b_{k+1} \| y_{k+1}^{ag}-\hat{x} \| ^2 -b_{k+1}(1-b_{k+1})\| y_{k}-y_{k+1}^{ag} \| ^2] \nonumber \\
            &\quad-(1-a_k-\theta_k^2 a_k-\delta)\|x_{k+1}-\bar{x}_{k}\|^2-(1-a_k)\|\bar{x}_{k}-x_{k}\|^2.   \label{l8}
        \end{align}
    Applying   $ \tau_{k}\theta_{k}\leq (1 +\theta_{k-1} )\tau_{k-1}$ to (\ref{l8}) means that
        \begin{align}
            &2\tau_k D_{\hat{x},\hat{y}}(y_{k+1}) + 2\tau_k(1+\theta_k)P_{\hat{x},\hat{y}}(x_{k+1}) \nonumber\\  
            &\quad +(1-a_k-\theta_k^2 a_k-\delta)\|x_{k+1}-\bar{x}_{k}\|^2+(1-a_k)\|\bar{x}_{k}-x_{k}\|^2 \nonumber \\
             &\quad+[a_{k+1} \| x_{k+1}^{ag}-x_{k+1} \| ^2 -a_{k+1} (1-a_{k+1})\| x_{k}-x_{k+1}^{ag} \| ^2] \nonumber\\
            &\quad +[\frac{1}{a_{k+2}}\| x_{k+2}^{ag}-\hat{x} \| ^2 -\frac{1-a_{k+2}}{a_{k+2}}\| x_{k+1}^{ag}-\hat{x} \| ^2+(1-a_{k+2})\| x_{k+1}-x_{k+1}^{ag} \|^2 ] \nonumber\\
            &\quad +\frac{1}{\beta}[(1-b_{k+1}-\delta) \| y_{k+1}-y_k \| ^2 +b_{k+1} \| y_{k+1}^{ag}-y_{k+1} \| ^2 -b_{k+1} (1-b_{k+1})\| y_{k}-y_{k+1}^{ag} \| ^2] \nonumber\\
            &\quad +\frac{1}{\beta}[\frac{1}{b_{k+2}}\| y_{k+2}^{ag}-\hat{x} \| ^2 -\frac{1-b_{k+2}}{b_{k+2}}\| y_{k+1}^{ag}-\hat{x} \| ^2+(1-b_{k+2})\| y_{k+1}-y_{k+1}^{ag} \|^2 ] \nonumber\\
            &\leq  2\tau_{k-1}(1+\theta_{k-1}) P_{\hat{x},\hat{y}}(x_{k-1}) \nonumber \\
            &\quad +[\frac{1-a_{k+1}}{a_{k+1}} \| x_{k+1}^{ag}-\hat{x} \| ^2 -\frac{(1-a_{k+1})^2}{a_{k+1}}\| x_{k}^{ag}-\hat{x} \| ^2+((1-a_{k+1})^2+a_k)\| x_{k}-x_{k}^{ag} \|^2 \nonumber \\
            &\quad + a_{k+1} \| x_{k+1}^{ag}-\hat{x} \| ^2 -a_{k+1}(1-a_{k+1})\| x_{k}-x_{k+1}^{ag} \| ^2]\nonumber \\
             &\quad +\frac{1}{\beta}[\frac{1-b_{k+1}}{b_{k+1}} \| y_{k+1}^{ag}-\hat{x} \| ^2 -\frac{(1-b_{k+1})^2}{b_{k+1}}\| y_{k}^{ag}-\hat{x} \| ^2+(1-b_{k+1})^2\| y_{k}-y_{k}^{ag} \|^2 \nonumber \\
            &\quad +b_{k+1} \| y_{k+1}^{ag}-\hat{x} \| ^2 -b_{k+1}(1-b_{k+1})\| y_{k}-y_{k+1}^{ag} \| ^2] . \label{l9}
        \end{align}
    Similarly to Theorem \ref{Th1}, from (\ref{l9}) we can deduce that $\{(x_k,y_k)\}$ are bounded sequences and
        \begin{align*}
            \lim_{k \rightarrow +\infty} \|y_k -y_{k-1} \| ^2 =0, \lim_{k \rightarrow +\infty} \|\bar{x}_k -x_{k} \| ^2 =0,\lim_{k \rightarrow +\infty}\| x_{k+1}^{ag}-x_{k+1} \|=0.
        \end{align*}
    Noticing that
         \begin{align*}
            \frac{x_{k+1}-x_{k}}{\tau_k}=\frac{\bar{x}_{k+1}-x_{k+1}}{\tau_{k+1}} \rightarrow 0 \quad as \; k \rightarrow +\infty,
        \end{align*}
    so we have 
        \begin{align*}
            \frac{x_{k+1}-x_{k+1}^{md}}{\tau_k}=[(1-a_{k+1})\frac{x_{k+1}-x_{k}}{\tau_k}+
            a_{k+1}\frac{x_{k+1}-x_{k+1}^{ag}}{\tau_k}] \rightarrow 0 \quad as \; k \rightarrow +\infty.
        \end{align*}
    Similarly,
        \begin{align*}
            \frac{y_{k+1}-y_{k+1}^{md}}{\beta \tau_k} \rightarrow 0 \quad as \; k \rightarrow +\infty.
        \end{align*}
    Let $\{(x_{k_i},y_{k_i})\}$ be a subsequence that converges to some cluster point $(x^*,y^*)$, we deduce $x_{k_i}^{ag}$ and $x_{k_i}^{md}$ converges $x^*$, $y_{k_i}^{ag}$ and $y_{k_i}^{md}$ converges $y^*$. Taking the limit in
        \begin{align*}
             \langle \frac{1}{\tau_{k_i}}(x_{{k_i}+1}-x_{{k_i}+1}^{md}) + K^T y_{{k_i}+1},x-x_{{k_i}+1}\rangle \; \geq \; g(x_{{k_i}+1})-g(x), \quad \forall x\in X ,\\
             \langle \frac{1}{\beta \tau_{k_i}} (y_{{k_i}+1}-y_{k_{i}+1}^{md})- K \bar{x}_{k_i} ,y-y_{{k_i}+1}\rangle \; \geq \; f^*(y_{{k_i}+1})-f^*(y), \quad \forall y\in Y,
        \end{align*}
    we obtain that $(x^*,y^*)$ is a saddle point of (\ref{p}).
    \par Setting $(\hat{x},\hat{y})=(x^*,y^*)$ in (\ref{l9}) and summing over from $k =1$ to $N$, similarly to Theorem \ref{Th1}    we have $x_N \rightarrow x^*$, $y_N \rightarrow y^*$ as $N \rightarrow +\infty$,
\end{proof}

\begin{thm}[ergodic convergence]
    Suppose that  $\{(x_k,y_k)\}$ is  generated by   NPDAL and $(\hat{x},\hat{y})$ is any saddle point of (\ref{p}).
    Then, it holds that
        \begin{align*}
            \mathcal{G}_{\hat{x},\hat{y}}(X_N,Y_N)=\mathcal{O}(\frac{1}{N}).
        \end{align*}
\end{thm}

\begin{proof}
   Summing up (\ref{l8}) from $k =1$ to $N$, we get 
        \begin{align}
            2\sum_{k=1}^{N}\omega_k \leq C_2,\label{l10}
        \end{align}
    where
        \begin{align*}
            C_2=&[(\frac{1}{a_2}+\frac{1}{a_3}+a_2-2-\frac{(1-a_3)^2}{a_3})\|x_2^{ag}-\hat{x}\|^2-\frac{(1-a_2)^2}{a_2}\|x_1^{ag}-\hat{x}\|^2] \nonumber\\
            &+\theta_1[(\frac{1}{a_1}+\frac{1}{a_2}+a_1-2-\frac{(1-a_2)^2}{a_2})\|x_1^{ag}-\hat{x}\|^2-\frac{(1-a_1)^2}{a_1}\|x_0^{ag}-\hat{x}\|^2] \nonumber\\
            &+\frac{1}{\beta}[(\frac{1}{b_2}+\frac{1}{b_3}+b_2-2-\frac{(1-b_3)^2}{b_3})\|y_2^{ag}-\hat{y}\|^2-\frac{(1-b_2)^2}{b_2}\|y_1^{ag}-\hat{y}\|^2].
        \end{align*}
    From the definition of $\omega_k$, we can know that left-hand side in (\ref{l10})  can be expressed as
        \begin{align}
            \sum_{k=1}^{N}\omega_k&=\tau_{N}(1+\theta_{N})P_{\hat{x},\hat{y}}(x_{N})
            +\sum_{k=2}^{N}[(1+\theta_{k-1})\tau_{k-1}-\tau_{k}\theta_{k}]P_{\hat{x},\hat{y}}(x_{k-1}) \nonumber\\
            & - \tau_{1}\theta_{1}P(x_{0})+\sum_{k=1}^{N}\tau_{k}D_{\hat{x},\hat{y}}(y_{k+1}).
    \end{align}
    By convexity of $P_{\hat{x},\hat{y}}(x)$ and $D_{\hat{x},\hat{y}}(y)$, we have    
        \begin{align}
            \tau_{N}(1+\theta_{N})P_{\hat{x},\hat{y}}(x_{N})&+\sum_{k=2}^{N}[(1+\theta_{k-1})\tau_{k-1}-\tau_{k}\theta_{k}]P_{\hat{x},\hat{y}}(x_{k-1}) \nonumber\\
            & \geq
            (\tau_{1}\theta_{1}+s_N)P_{\hat{x},\hat{y}}\left(\frac{\tau_{1}(1+\theta_{1} )x_1+\sum_{k=2}^{N}\tau_k \bar{x}_k}{\tau_{1}\theta_{1}+s_N}\right) \nonumber\\
            &=(\tau_{1}\theta_{1}+s_N)P_{\hat{x},\hat{y}}\left(\frac{\tau_{1}\theta_{1} 
            x_0+\sum_{k=1}^{N}\tau_k \bar{x}_k}{\tau_{1}\theta_{1}+s_N}\right)
            \geq
            s_N P_{\hat{x},\hat{y}}(X_{N}),
        \end{align}
    and
        \begin{align}
            \sum_{k=1}^{N}\tau_{k}D_{\hat{x},\hat{y}}(y_k)
            \geq s_N D_{\hat{x},\hat{y}}\left(\frac{\sum_{k=1}^{N}\tau_k y_k}{s_N}\right)
            =s_N D_{\hat{x},\hat{y}}(Y_N).
        \end{align}
    Hence,
        \begin{align*}
             \sum_{k=1}^{N}\omega_k\geq s_NP_{\hat{x},\hat{y}}(X_{N})+s_ND_{\hat{x},\hat{y}}(Y_N)-\tau_{1}\theta_{1}P_{\hat{x},\hat{y}}(x_{0}),
        \end{align*}
    and we have
        \begin{align*}
            \mathcal{G}_{\hat{x},\hat{y}}(X_N,Y_N)&=P_{\hat{x},\hat{y}}(X_{N})+D_{\hat{x},\hat{y}}(Y_N) \\
            &\leq  \frac{1}{ s_N}\left(C_2+\tau_{1}\theta_{1}P_{\hat{x},\hat{y}}(x_{0})\right).
        \end{align*}
   Together with $s_N=\sum_{k=1}^N\tau_k \geq \tau N$ implied by Lemma \ref{lemma3.1}, we obtain  $ \mathcal{G}_{\hat{x},\hat{y}}(X_N,Y_N)=\mathcal{O}(\frac{1}{N})$.  
\end{proof}

\section{Numerical experiment.}\label{section4}
  In the section, we conduct some numerical experiments to demonstrate the effectiveness  of three proposed algorithms NPDA, ANPDA and NPDAL by comparing them with some state-of-the-art  primal-dual algorithms on two groups of numerical examples  Matrix game and LASSO problems. The numerical experiments are performed by using Python 3 on   AMD Ryzen 5 2500U CPU 2.00GHz  and     Windows 10.   Codes can be found on https://github.com/Liuzexian2008/New-Primal-Dual-Algorithm-for-Convex-Problems.
\subsection{Matrix game.}
    We first conduct experiments on the following min-max matrix game problem:     
        \begin{equation} 
            \min_{x \in \bigtriangleup_n} \max_{y \in \bigtriangleup_m} \langle Kx,y\rangle,\label{N1}
        \end{equation}
    where  $K\in R^{m \times n}$ and $\bigtriangleup_p=\{z\in R^p: \sum_i z_i=1,z>0 \}$  is the standard unit simplex $in \; R^p$.
    \par Since neither $g$ nor $f^*$ is strongly convex for matrix game problem, we only study the performance of the PDA \cite{w1}, PDAL \cite{w2}, GRPDA \cite{w8}, GRPDAL \cite{w15}, NPDA and NPDAL.
    
    \par The initial points  are chosen as $x_0 = \frac{1}{n} (1,...,1)$ and $y_0 =\frac{1}{m} (1,...,1)$.  The matrix $K \in R^{m \times n}$ in     \eqref{N1} is generated  by the following ways \cite{w2} :  
        \begin{enumerate}[label=\arabic*.]
            \item $m = n = 100$. All entries of $K$ are generated independently from the uniform distribution in $[-1, 1]$.
            \item $m = n = 100$. All entries of $K$ are generated independently from the normal distribution $N (0, 1)$.
            \item $m = 500$, $n = 100$. All entries of $K$ are generated independently from the normal distribution $N (0, 1)$.
            \item $m = 1000$, $n = 2000$. The matrix $K$ is sparse with $10\%$
            generated independently from the uniform distribution in $[0, 1]$.
        \end{enumerate} 
Let   $\|K\|=\sqrt{\lambda_{\max}(K^TK)}$.  In the numerical experiment, the  parameters for each test algorithm are set as follows.
    \par PDA: $\tau = \sigma = 1/||K||$ \cite{w2}.
    \par PDAL: $\beta =1$, $\delta=0.99$, $\mu=0.7$, $\tau_0=\sqrt{\min \{m,n\}}/\|K\|_F$ \cite{w2}. 
    \par  GRPDA: $\psi=1.618$, $\tau = \sigma = \sqrt{\psi}/||K||$ \cite{w15}.
    \par GRPDAL: $\beta = 1$, $\sigma=0.99$, $\mu=0.7$, $\psi=1.5$, $\tau_0=\eta \sqrt\frac{\psi}{\beta}$, where $\eta=\frac{\|y_{-1}-y_0\| }{\|K^T(y_{-1}-y_0)\|}$\cite{w15}.
    \par NPDA: $\sigma = 1/||K||$, $\tau=\frac{2}{3}\sigma$, $a_k=b_k=0.01$.
    \par NPDAL: $\beta =1$, $\delta=0.96$, $\mu=0.7$, $\tau_0=\sqrt{\min \{m,n\}}/\|K\|_F$, $a_k=b_k=0.01$.

  The projection onto the unit simplex $\bigtriangleup_n$  or $\bigtriangleup_m$  is computed by the algorithm from \cite{w16}.         The primal-dual gap  
        \begin{align*} 
            \mathcal{G}(x,y)=\underset{i}{\max}(Kx)_i - \underset{j}{\min}(K^T y)_j  
        \end{align*}
        is used to measure the numerical performance of the test algorithms. 
    
       \begin{figure}[H] 
            \centering 
            \begin{minipage}{0.45\textwidth}
                \centering
                \includegraphics[width=\textwidth]{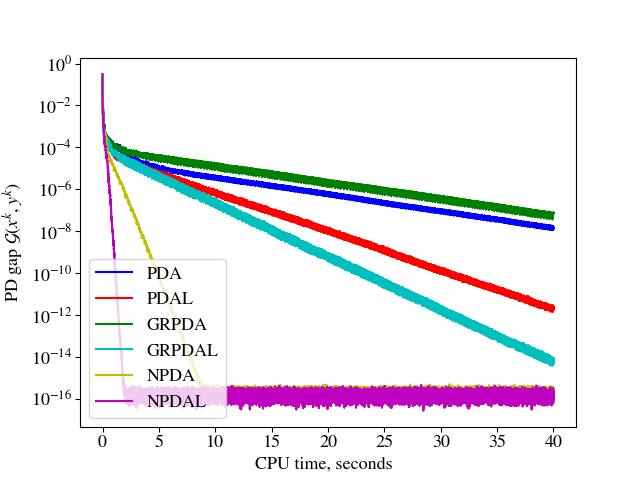}
                \subcaption{Example 1}
            \end{minipage}
        \hspace{0.05\textwidth}
            \begin{minipage}{0.45\textwidth}
            \centering
            \includegraphics[width=\textwidth]{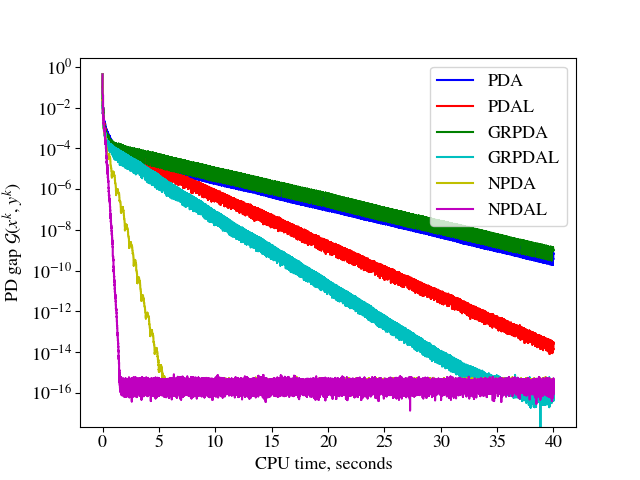}
            \subcaption{Example 2}
            \end{minipage}        
        \vspace{0.01\textheight} 
            \begin{minipage}{0.45\textwidth}
            \centering
            \includegraphics[width=\textwidth]{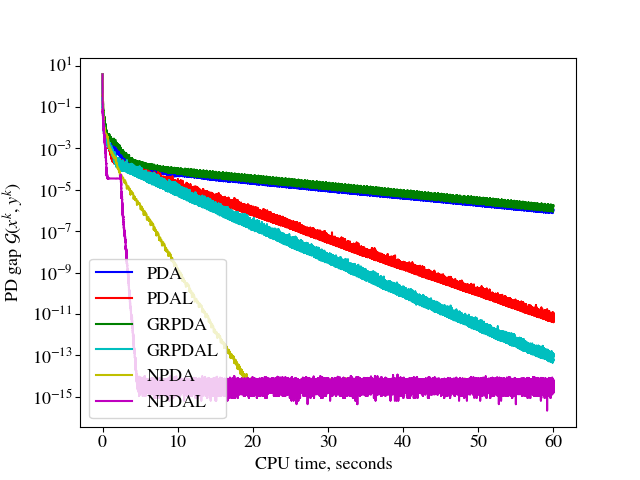}
            \subcaption{Example 3}
            \end{minipage}
        \hspace{0.05\textwidth}
            \begin{minipage}{0.45\textwidth}
            \centering
            \includegraphics[width=\textwidth]{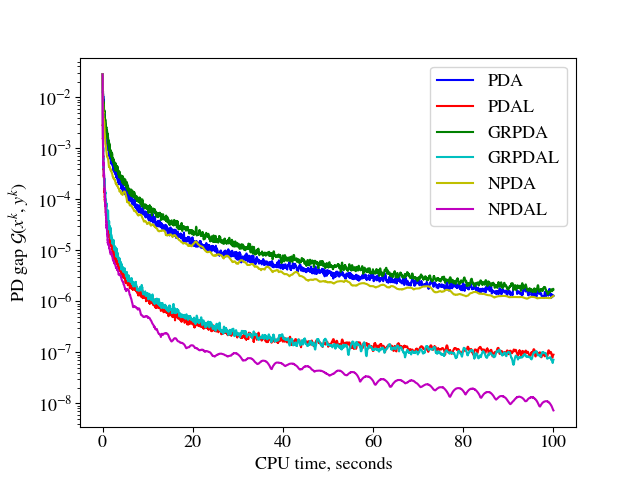}
            \subcaption{Example 4}
            \end{minipage}
        \caption{Matrix game}\label{Figure-1}
        \end{figure}
    \par The numerical results  for the primal-dual gap    $\mathcal{G}(x_k,y_k)$  at each iteration versus CPU time  are shown  in Figure \ref{Figure-1}. As shown in Figure \ref{Figure-1}, we observe that 
    NPDAL and NPDA have   significant advantages   over other test algorithms   for the first three examples, and NPDAL   performs better than   PDAL and GRPDAL for  all four example; GRPDAL is superior to  PDAL for the previous three examples  and is similar to  PDAL for the four example; PDA performs similarly or better than GRPDA.

\subsection{LASSO.}
We conduct the numerical experiment    on LASSO problem, which can be described as   the following $l_1$ regularized least squares problem:
        \begin{align}
            \underset{x \in R^{n}}{\min}\; \phi(x)=\frac{1}{2}\|Kx-b\|^2+\lambda\|x\|_1,\label{N2}
        \end{align}
    where $ K\in R^{m \times n}$ and  $b\in R^{m}$. Let $g(x)=\lambda\|x\|_1$, $f(z)=\frac{1}{2}\|Kz-b\|^2$. The problem    (\ref{N2}) can be written as the following minimax problem:
        \begin{align*}
            \underset{x\in R^{n}}{\min} \; \underset{y\in R^{m}}{\max}\;g(x)+\langle Kx,y \rangle-f^*(y),
        \end{align*}
    where $f^*(y)=\frac{1}{2}\|y+b\|^2-\frac{1}{2}\|b\|^2$.  It is not difficult to obtain that $\text{Prox}_{\lambda f^*(y)}=\frac{y-b\lambda}{1+\lambda}$.

   In the numerical experiment,   the matrix  $K\in R^{m \times n}$ in \eqref{N2} is constructed by   the following ways \cite{w2}:
        \begin{enumerate}
            \item $n = 1000$, $m = 200$, $s = 10$. All entries of $K$ are generated independently from $N(0, 1)$.
            \item $n = 2000$, $m = 1000$, $s = 100$. All entries of $K$ are generated independently from $N(0, 1)$.
            \item $n = 5000$, $m = 1000$, $s = 50$. First, we generate the matrix $B$ with entries from $N(0, 1)$. Then for any $p \in (0, 1)$ we construct the matrix $K$ by columns $K_j$, $j = 1, \ldots, n$, as follows: $K_1=\frac{B_1}{\sqrt{1-p^2}}, \; K_j=p*K_{j-1}+B_j$. As $p$ increases, $K$ becomes more ill-conditioned (where this example was considered). In this experiment we take $p = 0.5$.
            \item The same as the previous third example, but with $p = 0.9$.
        \end{enumerate}
 The vector  $b$ is set to $b = Kw + v$, where  $v \in R^m$ is generated by   $N (0, 0.1)$, and the  $s$   coordinates of   
        $w \in R^n$ are generated   by the uniform distribution in $[-10, 10]$ and the remaining  coordinates  are all zeros.     The initial points are $x_0 =  (0,...,0)$ and $y_0=K x_0-b$, and $\lambda=0.1$. 
        
        Let   $\|K\|=\sqrt{\lambda_{\max}(K^TK)}$.  
      The parameters for the test algorithms PDA, APDA, PDAL, GRPDA, AGRPDA,
GRPDAL, NPDA, ANPDA and NPDAL are set as follows. 
    \par PDA: $\sigma=\frac{1}{20\|K\|}$, $\tau=\frac{20}{\|K\|}$\cite{w2}.
    \par PDAL: $\beta=1/400$, $\delta=0.99$, $\mu=0.7$, $\tau_0=\sqrt{\min \{m,n\}}/\|K\|_F$\cite{w2}.
    \par APDA:  $\tau_0=\sigma_0=\frac{1}{\|K\|}$, $\gamma=0.1$\cite{w2}.
    \par GRPDA: $\psi=1.618$, $\sigma=\frac{\sqrt{\psi}}{20\|K\|}$, $\tau=\frac{20\sqrt{\psi}}{\|K\|}$\cite{w15}.
    \par AGRPDA: $\beta=1$,  $\tau=\frac{20\sqrt{\psi}}{\|K\|}$, $\gamma=0.01$\cite{w8}.
    \par GRPDAL: $\beta=400$, $\sigma=0.99$, $\mu=0.7$, $\psi=1.5$, $\tau_0=\eta \sqrt\frac{\psi}{\beta}$, where $\eta=\frac{\|y_{-1}-y_0\| }{\|K^T(y_{-1}-y_0)\|}$\cite{w15}.
    \par NPDA: $\sigma=\frac{3}{10\|K\|}$, $\tau=\frac{3}{\|K\|}$, $a_k=b_k=0.01$.
    \par ANPDA: $\sigma=\frac{3}{10\|K\|}$, $\tau=\frac{33}{10\|K\|}$, $a_k=b_k=0.005$.
    \par NPDAL: $\beta=1/10$, $\delta=0.96$, $\mu=0.7$, $\tau_0=\sqrt{\min \{m,n\}}/\|K\|_F$, $a_k=0.01$, $b_k=0$. \\       
   
\begin{figure}[H] 
            \centering 
            \begin{minipage}{0.45\textwidth}
                \centering
                \includegraphics[width=\textwidth]{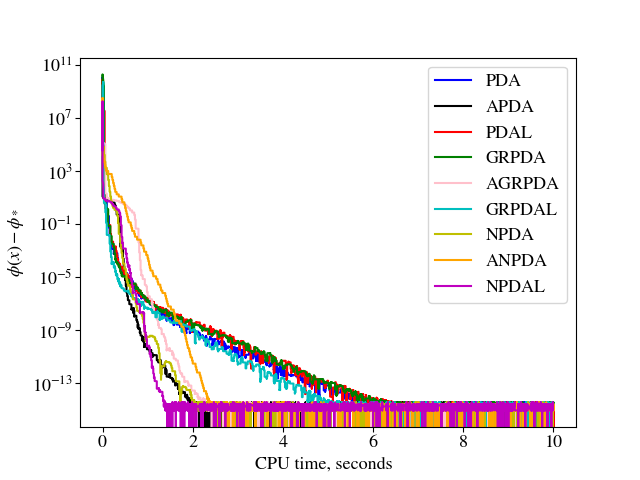}
                \subcaption{Example 1}
            \end{minipage}
        \hspace{0.05\textwidth}
            \begin{minipage}{0.45\textwidth}
            \centering
            \includegraphics[width=\textwidth]{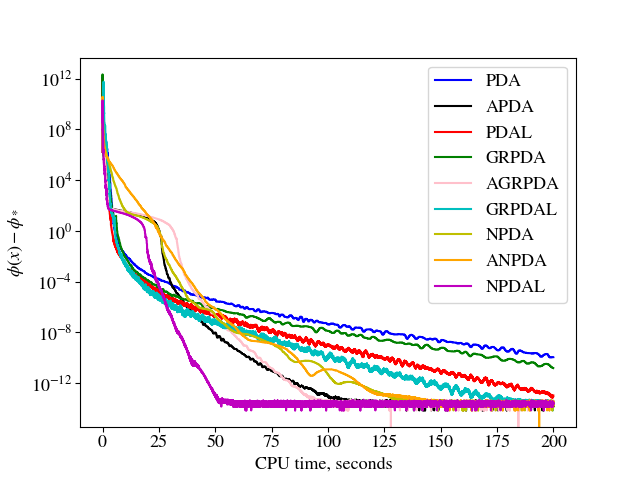}
            \subcaption{Example 2}
            \end{minipage}        
        \vspace{0.01\textheight} 
            \begin{minipage}{0.45\textwidth}
            \centering
            \includegraphics[width=\textwidth]{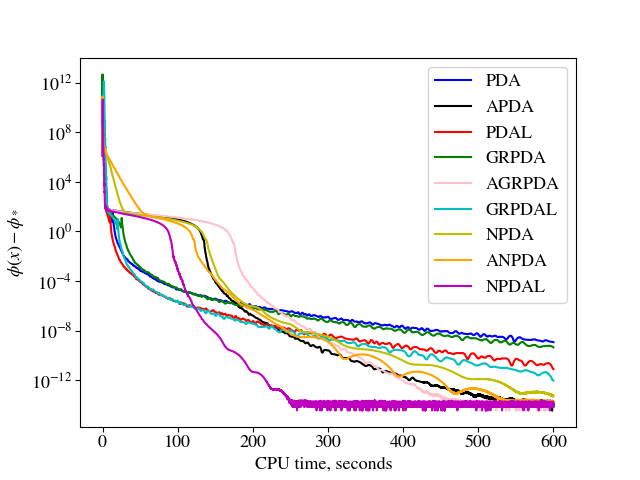}
            \subcaption{Example 3}
            \end{minipage}
        \hspace{0.05\textwidth}
            \begin{minipage}{0.45\textwidth}
            \centering
            \includegraphics[width=\textwidth]{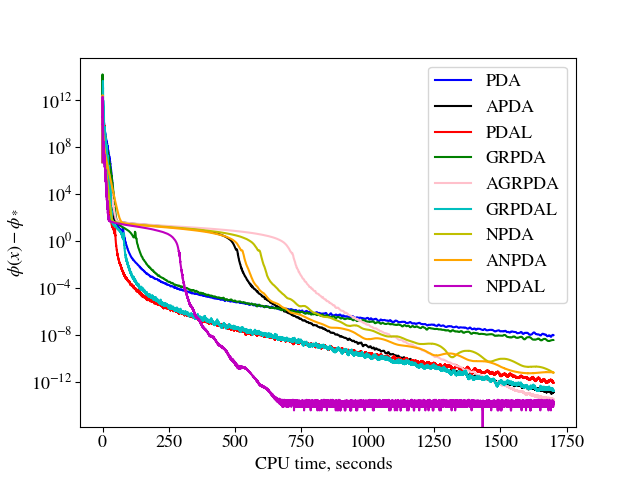}
            \subcaption{Example 4}
            \end{minipage}
        \caption{LASSO}\label{Figure-2}
        \end{figure}

   Since the   optimal solution of     (\ref{N2}) is unknown, we first run all the test algorithms for a sufficiently large number of iterations and then chose the minimum attainable function value  as an approximate solution $\phi_*=\phi(x_*)$ of (\ref{N2}). Figure \ref{Figure-2} collects the convergence results with $\phi(x_k)-\phi_*$  versus CPU time.   As shown in Figure \ref{Figure-2}, we can observe that NPDA always performs better than     PDA algorithm and   GRPDA algorithm for all numerical examples, and is even better than that of both PDAL and GRPDAL for the previous three numerical examples;  NPDAL achieves the best performance among all tested algorithms. 

\section{Conclusion.}\label{section5}
    Due to the significant importance of proximal terms for PDA,   we utilize the previous iterative points to exploit two new points $x_k^{md}$ and $y_k^{md}$  for   proximal terms and present a new  proximal-dual algorithm (NPDA) for convex-concave minimax problem with bilinear coupling term. The global convergence and convergence rate are also established.
     When either $f$ or $g$ is a strongly convex function, we  investigate an adaptive strategy for parameters for     NPDA and present  an  accelerated NPDA (ANPDA) with faster convergence rate $\mathcal{O}(1/N^2)$. Finally, we introduce  a linesearch strategy for  PDA  and establish its global convergence and convergence rate.
    \par Numerical experiments on matrix game problem and LASSO problem demonstrate that the two new proximal terms in NPDA can indeed bring large numerical improvements over PDA and  GRPDA,  and NPDAL performs  best among all the test algorithms for all numerical examples. \newline  \newline

 \noindent\textbf{Declaration of competing interest}\\
   The authors declare that they have no known competing financial interests or personal relationships that could have appeared to influence the work reported in this paper.\\

 \noindent\textbf{Acknowledgements}\\
 We   would like to thank Professor Yu-Hong Dai in Chinese
 Academy of Sciences for his valuable and insightful comments on this manuscript. This research is supported by the National Science Foundation of China (No. 12261019,12161053), Guizhou Provincial Science and Technology Projects (No. QHKJC-ZK[2022]YB084).\\

 \noindent\textbf{Data availability}\\
 The datasets generated during and analyzed during the current study are available from the corresponding author on reasonablerequest.\\

 \noindent\textbf{CRediT authorship contribution statement}\\

\textbf{Shuning Liu:} Writing - original draft, Investigation, Formal analysis, Conceptualization, Visualization. 
\textbf{Zexian Liu:} Writing - review \& editing, Conceptualization, Validation, Software, Methodology, Funding acquisition, Resources, Supervision, Validation.


\end{document}